\documentclass[11pt,a4paper]{article}
\usepackage{amsmath,amsfonts,amssymb,mathrsfs} 
\usepackage[english]{babel}
\newtheorem{theorem}{\bf Theorem}[section]
\newtheorem{lemma}[theorem]{\bf Lemma}

\newtheorem{definition}[theorem]{\bf Definition}
\newtheorem{remark}[theorem]{\bf Remark}
\newtheorem{proposition}[theorem]{ Proposition}

\newcommand{\abs}[1]{\lvert#1\rvert}
\newcommand{\norm}[1]{ \lVert {#1} \rVert}

\newcommand{\sob}[1]{L^{#1}(\mathbb{R}^n)}
\newcommand{\sobb}[2]{L^{#1}(#2)}
\newcommand{\sobolev}[1]{W^{1,#1}(\mathbb{R}^n)}

\newcommand{\integ}[2]{\int_{\mathbb{R}^n}  {#1} \, \mathrm{d}{#2}}
\newcommand{\integrale}[4]{\int_{#1}^{#2}  {#3} \, \mathrm{d}{#4}}
\newcommand{\scalare}[2]{\langle #1 , #2 \rangle}
\def \L {\mathscr{L}}
\def \LO {\mathscr{L}_0}
\def \K {K}
\def \R {{\mathbb {R}}}
\def \rn {{\mathbb {R}}^{N}}
\def \rnn {{\mathbb {R}}^{N+1}}
\def \div {{\text{\rm div}_X}}

\def\p{\partial}

\begin{document}
\title{A compactness result for the Sobolev embedding \\ via potential theory 
\footnote{AMS Subject Classification: 46E35, 35K70, 35B45, 35D30.}}
\author{{\sc{Filippo Camellini, Michela Eleuteri and Sergio Polidoro}}\\
Dipartimento di Scienze Fisiche Informatiche e Matematiche
\thanks{Modena via Campi 2013/b.
E-mail: michela.eleuteri@unimore.it, sergio.polidoro@unimore.it}
\thanks{Investigation supported by I.N.d.A.M.}
}

\maketitle

\footnotesize 
\begin{abstract}
In this note we give a proof of the Sobolev and Morrey embedding theorems based on the representation of functions 
in terms of the fundamental solution of suitable partial differential operators. We also prove the compactness of the 
Sobolev embedding. We first describe this method in the classical setting, where the fundamental solution of the 
Laplace equation is used, to recover the classical Sobolev and Morrey theorems. We next consider degenerate Kolmogorov 
equations. In this case, the fundamental solution is invariant with respect to a non-Euclidean translation group and 
the usual convolution is replaced by an operation that is defined in accordance with this geometry. We recover some 
known embedding results and we prove the compactness of the Sobolev embedding. We finally apply our regularity results 
to a kinetic equation.

\medskip\noindent
{\it Keywords: Sobolev spaces, Sobolev embedding, Morrey embedding, Compactness, Fundamental solution, Kolmogorov 
equation.}   
\end{abstract}

\normalsize

\section{Introduction}
\setcounter{section}{1} \setcounter{equation}{0}
\setcounter{theorem}{0} 

Sobolev and Morrey embedding theorems are fundamental tools in the regularity theory for Elliptic and Parabolic second 
order Partial Differential Equations (\emph{PDEs} in the sequel). In particular, they play a crucial role in the 
natural setting for the study of uniformly elliptic PDEs in divergence form, that is the Sobolev space $W^{1,p}$. 

There are several proofs of the Sobolev and Morrey embedding theorems, all of them rely on some integral representation 
of a general function $u \in W^{1,p}$ in terms of its gradient. Here we focus in particular on representation formulas 
based on the fundamental solution of the Laplace equation. 

Consider a function $u \in C_0^\infty(\mathbb{R}^n)$. By the very definition of fundamental solution $\Gamma$, the 
following identity holds
\begin{equation} \label{eq-fund-sol1}
  u(x) =  -  \int_{\mathbb{R}^n}\Gamma(x-y) \Delta u(y) \, d y, \qquad  \text{for every} \ x \in \mathbb{R}^n,
\end{equation}
and an integration by parts immediately gives
\begin{equation} \label{eq-fund-sol2}
  u(x) =   \int_{\R^n} \left\langle \nabla_y \Gamma(x-y), \nabla u(y) \right\rangle d y, \qquad 
  \text{for every} \ x \in \R^n,
\end{equation}
where $\langle \cdot, \cdot \rangle$ and $\nabla$ denote the usual inner product in $\R^n$ and the gradient, 
respectively. We recall that the gradient of the fundamental solution of the Laplace equation writes as follows 
\begin{equation} \label{eq-grad-fund-sol-lapl}
  \nabla \Gamma(x-y)= - \frac{1}{n \omega_n \abs{x-y}^n} (x-y), \qquad x \ne y,
\end{equation}
where $\omega_n$ is the measure of the $n$-dimensional unit ball. In particular, $\nabla \Gamma$ is an homogeneous 
function of degree $-n+1$, and there exists a positive constant $c_n$ such that
\begin{equation} \label{9Gamma}
\abs{\nabla \Gamma(x-y)}\leq c_n \abs{x-y}^{1-n}, 
\end{equation}
thus \eqref{eq-fund-sol2} yields the following inequality:
\begin{equation} \label{eq-conv-grad}
\abs{u(x)} \le c_n \int_{\R^n}{\abs{x-y}^{1-n}\abs{\nabla u(y)}} d {y}.
\end{equation}
The Young inequality for convolution with homogeneous kernels (see, for instance, Theorem 1, p. 119 in \cite{undici}) 
then gives
\begin{equation} \label{eq-Young}
  \| \nabla \Gamma * \nabla u \|_{L^{p^*}(\R^{n})} \le C_{p} \, \| \nabla u \|_{L^{p}(\R^{n})},  \qquad 1 < p < n,
\end{equation}
where $p^* = \frac{pn}{n-p}$ is the Sobolev conjugate of $p$, and $C_p$ is a positive constant which only depends on 
$p$ and on the dimension $n$. Here and in the sequel the dependence on $n$ will be often omitted. As a 
consequence we find
\begin{equation} \label{eq-pre-Sobolev}
  \| u \|_{L^{p^*}(\mathbb{R}^n)} \le C_{p} \, \| \nabla u \|_{L^{p}(\mathbb{R}^n)}, 
  \qquad \text{for every} \ u \in C_0^{\infty}(\mathbb{R}^n), \quad 1 < p < n.
\end{equation}
From the above inequality  we plainly obtain the following Sobolev inequality for any open set 
$\Omega \subseteq \R^{n}$ 
\begin{equation} \label{eq-Sobolev}
  \| u \|_{L^{q}(\Omega)} \le C_{p, q} \, \| u \|_{W^{1,p}(\Omega)}, 
  \qquad \text{for every} \ u \in W_0^{1,p}(\Omega), 
\end{equation}
with $1 < p < n$ and  $p \le q \le p^*$. Here $C_{p,q}$ is a positive constant which only depends on $p, q$ and $n$. 
By a standard argument \eqref{eq-Young} also gives the Sobolev embedding theorem for $W^{1,p}(\Omega)$ provided 
that the boundary of $\Omega$ is sufficiently smooth. 

The Morrey inequality (see Theorem \ref{e} below) can be obtained by the representation formula \eqref{eq-fund-sol2}, 
by using the following fact: there exists a positive constant $M_n$, only depending on $n$, such that 
%If $K$ is a $C^1$ function, defined on $\R^{n} \setminus \big\{ 0 \big\}$ and is homogeneous of degree $- \alpha$, 
with 
%$0 < \alpha < n$, then there exists a positive constant $M$ such that
\begin{equation}\label{eq-grad-hom}
  | \partial_{x_j} \Gamma (x) - \partial_{x_j} \Gamma (y)  | \le M_n \frac{|x-y|}{|x|^{n}}, 
  \quad \text{for} \ j= 1, \dots, n,
\end{equation}
for every $x,y \in \R^{n} \setminus \big\{ 0 \big\}$ such that $|x-y| \le |x|/2$. Indeed, a rather simple argument 
based on \eqref{eq-grad-hom} provides us with the following bound: if $u \in W_0^{1,p}(\Omega)$, 
with $p > n$, then 
\begin{equation} \label{eq-Morrey}
  | u(x) - u(y) | \le \widetilde C_{p} \, \| \nabla u \|_{L^{p}( \Omega)} |x-y|^{1 - \frac{n}{p}}, 
  \qquad \text{for every} \ x, y \in\Omega,
\end{equation}
for some positive constant $\widetilde C_{p}$ only depending on $p$ and $n$.

It is worth noting that the inequality \eqref{eq-grad-hom} can be also used to prove the compactness of the Sobolev 
embedding \eqref{eq-Sobolev} for $p < q < p^*$, if $\Omega$ is a bounded open set. As we will see in the sequel, the 
following estimates holds for $p < q < p^*$: there exists a positive constant $\widetilde C_{p, q}$ such that 
\begin{equation} \label{eq-Sobolev-unif}
  \| u ( h + \cdot) - u \|_{L^{q}(\Omega)} \le \widetilde C_{p, q} \, \| \nabla u \|_{L^{p}(\Omega)} \
  |h|^{n \left( \frac{1}{q}- \frac{1}{p^*} \right)}, 
\end{equation}
for every $u \in C_0^{\infty}(\Omega)$ and for every $h \in \R^n$ sufficiently small. Note that the exponent in the 
right hand side of \eqref{eq-Sobolev-unif} belongs to the interval $]0,1[$ if, and only if, $p < q < p^*$, 
then in this case we have %is strictly positive
\begin{equation*} %\label{eq-Sobolev-unif-2}
  \| u ( h + \cdot) - u \|_{L^{q}(\Omega)} \to 0 \qquad \text{as} \qquad |h| \to 0.
\end{equation*}
This inequality provides us with the \emph{integral uniform continuity}, which  is needed for the compactness in the 
$L^q$ spaces. We also observe that $n \Big( \dfrac{1}{q}- \dfrac{1}{p^*} \Big) \to 1$ as $q \rightarrow p$. We 
then retrieve a known result contained for instance in \cite{undici}, Chapter V, Section 3.5.

\medskip

The advantage of the method described above, with respect to other ones, is in that it only requires the 
existence of a \emph{fundamental solution} and its \emph{homogeneity} properties. In particular, it applies to the 
function spaces introduced by Folland \cite{nove} for the study second order linear differential operators that satisfy 
the H\"ormander's condition (see \cite{dieci}). It should be noticed that this approach has also a drawback, in  
that it does not provide us with the Sobolev inequality for $p=1$. On the other hand it is unifying, as it gives 
the Sobolev and Morrey embedding theorems and a compactness result by using a single representation formula. 

To clarify the use of this method to the study of the so-called H\"ormander's operators we next focus on the degenerate 
Kolmogorov $\L_0$ on $\R^{2n+1}$, which is one of the simplest examples belonging to this class. Let $\Omega$ be an 
open 
subset of $\R^{2n+1}$ and let $u$ be a smooth real valued function defined on $\Omega$. We denote the variable of 
$\R^{2n+1}$ as follows $z = (x,y,t) \in \R^n \times  \R^n \times  \R$, and we set
\begin{equation}\label{e-Kolmo-0}
  \L_0 u := \Delta_{x} u + \langle x, \nabla_y u \rangle - \partial_t u, \qquad 
  \Delta_{x} u := \sum_{j=1}^n \partial_{x_j}^2 u.
\end{equation}
As we will see in the sequel (see equation \eqref{eq-fund-sol-K} below) the function $\Gamma$ defined as
\begin{equation*}
\begin{cases}
\displaystyle
\Gamma(x,y,t) = \frac{\widetilde c_n}{t^{2n}} \exp\left( - \tfrac{|x|^2}{t} - 3 \tfrac{\langle x, y 
\rangle}{t^2} - 3 \tfrac{|y|^2}{t^3}  \right), &  \text{for} \ (x,y,t) \in \R^{2n} \times ]0, + \infty[, \\
\Gamma(x,y,t) = 0, & \text{for} \ (x,y,t) \in \R^{2n} \times ]- \infty, 0],
\end{cases}
\end{equation*}
is the fundamental solution of $\L_0$. Here $\widetilde c_n = \frac{{3}^{n/2}}{(2 \pi)^n}$. In particular, in analogy 
with the heat equation, we have that the function $u$ defined as 
\begin{equation} \label{eq-fund-sol-K0}
\begin{split}
 u(x,y,t) = & \int_{\R^{2n}} \Gamma(x - \xi,y + t \xi -\eta,t- t_0) \varphi(\xi, \eta) d \xi \, d \eta \, - \\
 & \int_{\R^{2n} \times ]t_0,t[}  \Gamma(x - \xi,y + (t- \tau) \xi -\eta,t- \tau) 
 f(\xi, \eta, \tau) d \xi \, d \eta \, d \tau
\end{split}
\end{equation}
is a solution to the following Cauchy problem
\begin{equation*}
\begin{cases}
 \L_0 u = f  & \text{in} \quad \R^{2n} \times ]t_0, + \infty[, \\
 u_{\mid t=t_0} = \varphi & \text{in} \quad \R^{2n}.
\end{cases}
\end{equation*}
whenever $f$ and $\varphi$ are bounded continuous functions. 

A remarkable fact is that a kind of convolution is hidden in the expression \eqref{eq-fund-sol-K0}. More specifically, 
we define the operation ``$\circ$'' by setting
\begin{equation} \label{eq-group-K0}
  (x,y,t) \circ (\xi, \eta, \tau) := (x + \xi, y + \eta + \tau x, t + \tau), \qquad (x,y,t), (\xi, \eta, \tau) \in 
\R^{2n+1},
\end{equation}
and we note that $\left( \R^{2n+1}, \circ \right)$ is a non commutative group. The identity of the group is $(0,0,0)$ 
and the inverse of $(x,y,t)$ is $(-x, -y + xt, -t)$. With this notation, it is easy to check that the expression 
appearing in \eqref{eq-fund-sol-K0} can be written as follows
\begin{equation*}
  (x - \xi,y + (t- \tau) \xi -\eta,t- \tau) =  (\xi, \eta, \tau)^{-1} \circ (x,y,t).
\end{equation*}
Moreover, the group $\left( \R^{2n+1}, \circ \right)$ is homogeneous with respect to the dilation defined as $d_r 
(x,y,t) := \left( r x, r^3 y, r^2 t\right)$, in the sense that 
\begin{equation} \label{eq-dil-K0}
  d_r \big( (x,y,t) \circ (\xi, \eta, \tau) \big) = d_r (x,y,t) \circ d_r (\xi, \eta, \tau), 
  \qquad (x,y,t), (\xi, \eta, \tau) \in \R^{2n+1}, r>0.
\end{equation}
This algebraic structure was introduced and studied by Lanconelli and Polidoro in \cite{dodici}. In \cite{dodici} 
it was also noticed that $\Gamma$ is homogeneous of degree $- 4n$ with respect to $\left( d_t \right)_{r>0}$, that is
\begin{equation} \label{eq-dilGamma-K0}
  \Gamma \big(d_r (x,y,t) \big) = \frac{1}{r^{4n}} \Gamma (x,y,t), \qquad (x,y,t) \in \R^{2n+1}, r>0.
\end{equation}
Moreover, if we let $z = (x,y,t), \zeta = (\xi, \eta, \tau)$, then \eqref{eq-fund-sol-K0} can be written as follows
\begin{equation} \label{eq-fund-sol-K0-bis}
 u(z) = \int_{\R^{2n}} \Gamma((\xi,\eta,t_0)^{-1} \circ z) \varphi(\xi, \eta) d \xi \, d \eta \, - \int_{\R^{2n} 
\times ]t_0,t[}  \Gamma(\zeta^{-1} \circ z)  f(\zeta) d \zeta.
\end{equation}
In particular, if $u \in C_0^\infty(\R^{2n+1})$ and supp$(u) \subset \big\{ t > t_0 \big\}$, then we have that
\begin{equation} \label{eq-fund-sol1-bis}
  u(z) =  -  \int_{\mathbb{R}^{2n+1}}\Gamma(\zeta^{-1} \circ z) \L_0 u(\zeta) \, 
  d \zeta, \quad \text{for every} \ z \in \R^{2n+1},
\end{equation}
which is analogous to \eqref{eq-fund-sol1}. Summarizing: the operation in \eqref{eq-fund-sol1-bis} is considered here 
as 
a \emph{convolution} with respect to the non-Euclidean operation ``$\circ$'' defined in \eqref{eq-group-K0}, with a 
kernel $\Gamma$ that is homogeneous whit respect to the \emph{anisotropic} dilation $d_r$. Based on this representation 
formula, we prove Sobolev and Morrey theorems for solutions to Kolmogorov equations in divergence form $\L u = 
\text{\rm 
div}_{x} F + f$, where
\begin{equation}\label{e-Kolmo-1}
  \L u := \text{\rm div}_{x} \left(A(z) \nabla_x u \right) + \langle x, \nabla_y u \rangle - \partial_t u.
\end{equation}
Here $A$ is a $n \times n$ symmetric matrix with bounded and measurable coefficients and, for every vector field $F 
\in C^1(\R^{2n+1}, \R^n)$ we denote $\text{\rm div}_{x} F(x,y,t) := \sum_{j=1}^n \partial_{x_j} F_j(x,y,t)$. 
In order to simplify our treatment, we suppose that $F = 0$ and $f=0$, so that $u$ is a solution of $\L u = 0$. 
In this case we have that $\L_0 u = \text{\rm div}_{x}(I_n - A) \nabla_x u$, where $I_n$ denotes the $n \times n$ 
identity matrix. Then, an integration by parts in \eqref{eq-fund-sol1-bis} gives
\begin{equation} \label{eq-fund-sol2-bis}
  u(z) =   \int_{\mathbb{R}^{2n+1}} \langle (I_n - A(\zeta)) \nabla_\xi \Gamma(\zeta^{-1} \circ z),
  \nabla_{\xi} u(\zeta)\rangle \, d \zeta,
\end{equation}
for every solution $u$ to $\L u = 0$. It is known that the derivatives $\partial_{\xi_1} \Gamma, \dots, 
\partial_{\xi_n} \Gamma$ are homogeneous functions of degree $- (2n+1)$ with respect to the dilation $\left( d_r 
\right)_{r>0}$. Moreover, the coefficients of the matrix $I_n - A$ are bounded, then the above identity provides us 
with 
the analogous of \eqref{eq-fund-sol2} for the solutions $u$ to the equation $\L u = 0$. 

We point out that only the derivatives with respect to the first $n$ variables of the gradient of $u$ appear in 
the representation formula \eqref{eq-fund-sol2-bis}, then a Sobolev inequality holding for all functions cannot 
be obtained from \eqref{eq-fund-sol2-bis}, because of the lack of information on the remaining $n$ direction. 
Nevertheless, this formula is used by Cinti, Pascucci and Polidoro in  \cite{sei,sette} to prove a Sobolev embedding 
theorem for \emph{solutions} to the Kolmogorov equation $\L u = 0$.  Indeed, in \cite{sei,sette} the Sobolev theorem 
for 
\emph{solutions} is combined  with a Caccioppoli inequality, still for \emph{solutions}, in order to apply the Moser's 
iterative method and prove an $L^\infty_{{\rm loc}}$ estimate for the solutions to $\L u = 0$. 
We also recall that a Morrey result for the solutions to $\L u = \text{\rm div}_{x} F$ was proven by Manfredini and 
Polidoro in \cite{otto}, and later by Polidoro and Ragusa in \cite{otto.1} by the same method used here.  

In this note we are concerned with the compactness of the Sobolev embedding for the solutions to $\L u = 0$ for a 
family of degenerate Kolmogorov equations, defined on $\mathbb{R}^{N+1}$, that will be still denoted by $\L$. As 
we will see in Section 3, the operator \eqref{e-Kolmo-1} is the prototype of this family of degenerate operators, and 
in this case, $N = 2n$. In Section 3 we introduce the notation that will be used in the following part of this 
introduction, and we will state the conditions (H.1) and (H.2) that ensure that the \emph{principal part} $\L_0$ of 
$\L$ has a smooth fundamental solution $\Gamma$, which is invariant with respect to a translation analogous to 
\eqref{eq-group-K0}, and homogeneous of degree $-Q$, with respect to a dilation analogous to \eqref{eq-dil-K0}. We will 
refer to the positive integer $Q+2$ as \emph{homogeneous dimension} of the space $\rnn$ and plays the role of $n$ in 
the 
Euclidean setting $\R^n$ where the elliptic operators are studied. In the sequel $p^*$ and $p^{**}$ denote the positive 
numbers such that 
\begin{equation} \label{eq-p-pstar}
  \frac{1}{p^*} = \frac{1}{p} - \frac{1}{Q+2}, \qquad \frac{1}{p^{**}} = \frac{1}{p} - \frac{2}{Q+2}.
\end{equation}
Clearly, $p^*$ and $p^{**}$ are finite and positive whenever $1 \le p < Q+2$ and $1 \le p < \dfrac{Q+2}{2}$, 
respectively. 

Our main result is the following Theorem. It provides us with some estimates of the convolution of a function 
belonging to some $L^p$ space with the fundamental solution $\Gamma$ and with its derivatives $\partial_{x_j} \Gamma$, 
$j = 1, \dots, m_0$, with $m_0 \le N$. These estimates, applied to the representation formula for solutions to 
$\L u = 0$ given in Theorem \ref{th-rappr}, yield Sobolev theorems, Morrey theorems and the compactness of the Sobolev 
embedding.

\begin{theorem} \label{main-th-1}
Let $\L$ be an operator in the form \eqref{ultraparabolica}, satisfying the hypotheses {\rm (H.1)} and {\rm (H.2)} in 
Section 3, and let $\Gamma$ be the fundamental solution of its \emph{principal part}. Let also $Q+2$ be the homogeneous 
dimension of the space $\mathbb{R}^{N+1}$, and let $p$ be such that $1 \le p < + \infty$. For every $f, g_j \in 
L^p(\rnn)$ we let $u, v_j$ be defined as follows
\begin{equation*} %\label{eq-convolution-main}
  u(z) = \int_{\rnn} \Gamma(\zeta^{-1} \circ z) f(\zeta) \, {\rm d} \zeta,
  \quad v_j(z) = \int_{\rnn} \partial_{x_j} \Gamma(\zeta^{-1} \circ z) g_j (\zeta) \, {\rm d} \zeta,
  \quad j = 1, \dots, m_0.
\end{equation*}
Then, for every $j = 1, \dots, m_0$ we have: 
\begin{itemize}
  \item (Sobolev) if $1 < p < Q+2$, then there exists a positive constant $C_p$ such that
\begin{equation*} %\label{eq-Young-main}
  \| v_j \|_{L^{p^*}(\rnn)} \le C_{p} \, \|g_j \|_{L^{p}(\rnn)}, 
\end{equation*}
 \item (Compactness) if moreover $p < q < p^*$, then there exists a positive constant 
$\widetilde C_{p, q}$ such that 
\begin{equation*} %\label{eq-Sobolev-unif-main}
  \| v_j ( \cdot \circ h) - v_j \|_{L^{q}(\rnn)} \le \widetilde C_{p, q} \, \| g_j \|_{L^{p}(\rnn)} \
  \|h\|^{(Q+2) \left( \frac{1}{q}- \frac{1}{p^*} \right)}, 
\end{equation*}
for every $h \in \rnn$,
  \item (Morrey) if  $p > Q+2$, then there exists a positive constant $\widetilde C_p$ such that
\begin{equation*} %\label{eq-Morrey-main}
   | v_j(z) - v_j(\zeta) | \le \widetilde C_{p} \, \|  g_j \|_{L^{p}( \rnn)} \|\zeta^{-1} \circ z\|^{1 - 
\frac{Q+2}{p}}, \quad \text{for every} \ z, \zeta \in \rnn.
\end{equation*}
 \end{itemize}
 We also have
\begin{itemize}
  \item (Sobolev) if $1 < p < \dfrac{Q+2}{2}$, then there exists a positive constant $C_p$ such 
that
\begin{equation*} %\label{eq-Young-main}
  \| u \|_{L^{p^{**}}(\rnn)} \le C_{p} \, \| f \|_{L^{p}(\rnn)}, 
\end{equation*}
  \item (Compactness) if $p^* < q < p^{**}$, then there exists a positive constant $\widetilde 
C_{p, q}$ such that 
\begin{equation*} %\label{eq-Sobolev-unif-main}
  \| u ( \cdot \circ h) - u \|_{L^{q}(\rnn)} \le \widetilde C_{p, q} \, \| f \|_{L^{p}(\rnn)} \
  \|h\|^{(Q+2) \left( \frac{1}{q}- \frac{1}{p^{**}} \right)}, 
\end{equation*}
for every $h \in \rnn$,
  \item (Morrey) if $\dfrac{Q+2}{2} < p < Q+2$, then there exists a positive constant $\widetilde C_p$ 
such that
\begin{equation*}
     | u(z) - u(\zeta) | \le \widetilde C_{p} \, \|  f \|_{L^{p}( \rnn)} \|\zeta^{-1} \circ z\|^{2 - 
\frac{Q+2}{p}},   \qquad \text{for every} \ z, \zeta \in \rnn.
\end{equation*}  
\end{itemize}
\end{theorem}

From the above result and a representation formula for the solution to $\L u = 0$ we obtain the following result. 

\begin{theorem} \label{main-th-2} Let $\Omega$ be an open set of $\rnn$, and let $u$ be a weak solution to  
$\L u = 0$ in $\Omega$. Suppose that $u, \partial_{x_1} u, \dots, \partial_{x_{m_0}} u \in L^p%_\text{\rm loc}
(\Omega)$. Then for every compact set $\K \subset \Omega$, there exist a positive constant $\widetilde \varrho$ 
such that we have:
\begin{itemize}
  \item (Sobolev embedding) if $1 < p < Q+2$, then there exists a positive constant $C_p$ such that
\begin{equation*} %\label{eq-Young-main}
   \| u \|_{L^{p^*}(\K)} \le \,  C_{p} \bigg( \| u \|_{L^{p}(\Omega)} + \sum_{j=1}^{m_0}\| \partial_{x_j} u 
\|_{L^{p}(\Omega)} \bigg),
\end{equation*}
 \item (Compactness) if moreover $p < q < p^*$, then there exists a positive constant $\widetilde C_{p, q}$ such that 
\begin{equation*} %\label{eq-Sobolev-unif-main}
  \| u ( \cdot \circ h) - u \|_{L^{q}(\K)} \le \widetilde C_{p, q}  
  \bigg( \| u \|_{L^{p}(\Omega)} + \sum_{j=1}^{m_0}\| \partial_{x_j} u \|_{L^{p}(\Omega)} \bigg) \
  \|h\|^{(Q+2) \left( \frac{1}{q}- \frac{1}{p^*} \right)},
\end{equation*}
for every $ h \in \rnn$ such that $\| h \| \le \widetilde \varrho$, 
\item (Morrey embedding) if  $p > Q+2$, then there exists a positive constant $\widetilde C_p$ such that
\begin{equation*} %\label{eq-Morrey-main}
   | u (z) - u (\zeta) | \le  \widetilde C_{p}  
   \bigg( \| u \|_{L^{p}(\Omega)} + \sum_{j=1}^{m_0}\| \partial_{x_j} u \|_{L^{p}(\Omega)} \bigg)
   \|\zeta^{-1} \circ z\|^{1 - \frac{Q+2}{p}},
\end{equation*}
for every $\ z, \zeta \in \K$ such that $\|\zeta^{-1} \circ z\| \le \widetilde \varrho$.
 \end{itemize}
\end{theorem}

\medskip

The following Theorem is related to the main result of the article \cite{quattordici} by Bouchut, where the regularity 
of the solution of the kinetic equation
\begin{equation} \label{eq-kinetic}
  \partial_t f + \langle v, \nabla_x f \rangle = g, \qquad (t,x,v) \in \Omega \subseteq \R \times \R^{n} \times \R^{n},
\end{equation}
is considered. Note that the differential operator appearing in the left hand side of \eqref{eq-kinetic} agrees 
with the first order part of $\L$ defined in \eqref{e-Kolmo-1}. Actually, the notation of the following result refers 
to this operator, and, in particular, the homogeneous dimension of the space $\R^{2n + 1}$ is in this case $Q+2 = 4 
n + 2$.
%domain of the solution $f$ is any open set $\Omega \subset \R^{2n + 1}$. Moreover,  

\medskip

\begin{theorem} \label{main-th-3} Let $\Omega$ be an open set of $\R^{2n+1}$, and let $f \in L^2_\text{\rm 
loc}(\Omega)$ be a weak solution to  \eqref{eq-kinetic}. Suppose that $g, f, \partial_{v_1} f, \dots, \partial_{v_{n}} 
f \in L^p%_\text{\rm loc}
(\Omega)$. Then for every compact set $\K \subset \Omega$, there exist a positive constant 
$\widetilde \varrho$ 
such that we have:
\begin{itemize}
  \item (Sobolev embedding) if $1 < p < 4n+2$, then there exists a positive constant $C_p$ such that
\begin{equation*} %\label{eq-Young-main}
   \| f \|_{L^{p^*}(\K)} \le \,  C_{p} \bigg( \| g \|_{L^{p}(\Omega)} + \| f \|_{L^{p}(\Omega)} + \sum_{j=1}^{n}\| 
\partial_{v_j} f \|_{L^{p}(\Omega)} \bigg),
\end{equation*}
 \item (Compactness) if moreover $p < q < p^*$, then there exists a positive constant $\widetilde C_{p, q}$ such that 
\begin{equation*} %\label{eq-Sobolev-unif-main}
  \| f ( \cdot \circ h) - f \|_{L^{q}(\K)} \le \widetilde C_{p, q}  
  \bigg( \| g \|_{L^{p}(\Omega)} + \| f \|_{L^{p}(\Omega)} + \sum_{j=1}^{n}\| \partial_{v_j} f \|_{L^{p}(\Omega)} 
\bigg) 
  \ \|h\|^{(4n + 2) \left( \frac{1}{q}- \frac{1}{p^*} \right)},
\end{equation*}
for every $ h \in \R^{2n+1}$  such that $\| h \| \le \widetilde \varrho$, 
\item (Morrey embedding) if  $p > 4n + 2$, then there exists a positive constant $\widetilde C_p$ such that
\begin{equation*} %\label{eq-Morrey-main}
   | f (z) - f (\zeta) | \le  \widetilde C_{p}  
   \bigg( \| g \|_{L^{p}(\Omega)} + \| f \|_{L^{p}(\Omega)} + \sum_{j=1}^{n}\| \partial_{v_j} f \|_{L^{p}(\Omega)} 
   \bigg)  \|\zeta^{-1} \circ z\|^{1 - \frac{4n + 2}{p}},
\end{equation*}
for every $\ z, \zeta \in \K$ such that $\|\zeta^{-1} \circ z\| \le \widetilde \varrho$.
 \end{itemize}
\end{theorem}

The proof of Theorems \ref{main-th-2} and \ref{main-th-3} is given in Section 4. 

\medskip

We next give some comments to our main results. We still refer here to the notation relevant to the operator $\L$ 
defined in \eqref{e-Kolmo-1}, and to the representation formula \eqref{eq-fund-sol2-bis}. As we said above, it 
holds for solutions to $\L u = 0$ then, for this reason, it seems to be weaker than the usual Sobolev inequality. On 
the other hand, due to the strong degeneracy of the operator $\L$, its natural Sobolev space $W^{1,p}_{\L}$ is the space 
of the functions $u \in L^p$ with weak derivatives $\partial_{x_1} u, \dots, \partial_{x_n} u \in L^p$. In particular, 
it is impossible to prove a Sobolev inequality unless some information is given on $u$ with respect to the remaining 
variables $y_1, \dots, y_n$ and $t$. We obtain this missing information from the fact that $u$ is a solution to $\L u = 
0$ (or, in a more general case, to $\L u = \text{\rm div}_{x} F + f$). We also note that the regularity property of the 
operator $\L$ is quite unstable. Indeed, let us fix any $x_0 \in \R^n$ and consider the operator $\widetilde \L_0$, 
acting on $(x,y,t) \in \R^{2n+1}$ as follows
\begin{equation*} %\label{e-Kolmo-01}
  \widetilde \L_0 u := \Delta_{x} u + \langle x_0, \nabla_y u \rangle - \partial_t u.
\end{equation*}
Its natural Sobolev spaces agrees with that of $\L$, however it is known that a fundamental solution for $\widetilde 
\L_0$ does not exists and our method for the proof of the Sobolev inequality fails in this case. Actually, it is not 
difficult to check that the Sobolev inequality does not hold for the solutions to $\widetilde \L_0 u = 0$.

We conclude this discussion with a simple remark. Also when we consider the more familiar uniformly parabolic 
equations, we find that the natural Sobolev space only contains the \emph{spatial} derivatives, and it is not possible 
to find a simple natural space for the time derivative. As a matter of facts, several regularity results for parabolic 
equations depend on some \emph{fractional Sobolev spaces}. The situation becomes more complicated when we consider 
second order PDEs with non-negative characteristic form analogous to $\L$. An alternative approach to our method, that 
only relies on a representation formula in terms of the fundamental solution, is the use of \emph{fractional Sobolev 
spaces} (we refer to the articles by Bochut \cite{quattordici}, see also Golse, Imbert, Mouhot and Vasseur 
\cite{tredici}) to recover the missing information with respect to the variables $y_1, \dots, y_n$ and $t$.

\medskip

This article is organized as follows. In Section 2 we give a comprehensive proof of the Sobolev embedding, of its 
compactness, and the Morrey embedding, following the method above outlined. In Section 3 we recall the tools of the 
Real Analysis on Lie groups we need to prove Theorem \ref{main-th-1}, and we give its proof. In Section 4 we discuss 
some  applications of Theorem \ref{main-th-1} to the solutions of $\L u = 0$. Section 5 contains some comments about 
the possible extension of Theorem \ref{main-th-1} to a family of more general operators considered by Cinti and 
Polidoro in \cite{quindici}. 

\section{Continuous and compact embeddings: the Euclidean case}
\setcounter{section}{2} \setcounter{equation}{0}
\setcounter{theorem}{0} 

In this Section we give a comprehensive proof of the Sobolev embedding \eqref{eq-Sobolev}, the Morrey embedding
\eqref{eq-Morrey}, and of the inequality \eqref{eq-Sobolev-unif} from which the compactness of the Sobolev 
embedding follows. As said in the Introduction, all these results rely on the representation formula 
\eqref{eq-fund-sol2}, which gives the bound \eqref{eq-conv-grad} that we recall below
\begin{equation*} %\label{eq-conv-grad}
  \abs{u(x)} \le c_n \int_{\R^n}{\abs{x-y}^{1-n}\abs{\nabla u(y)}} \, d {y}, 
  \qquad \text{for every} \ u \in C_0^\infty (\R^n).
\end{equation*}
With this aim, we first recall the weak Young inequality that gives the Sobolev embedding, then we prove 
\eqref{eq-grad-hom} and we deduce from this and \eqref{eq-conv-grad} the Morrey embedding
\eqref{eq-Morrey}, and that stated in the inequality \eqref{eq-Sobolev-unif}.

\subsection{Some preliminary results} 

For a given positive $\alpha$ %such that $0 < \alpha < n$ 
we denote by $K_\alpha$ any continuous homogeneous function of degree $- \alpha$, that is a function satisfying
\begin{equation*}
  K_\alpha(r x)= r^{- \alpha } K_\alpha (x), \quad 
  \text{for every} \  x \in \mathbb{R}^n \setminus \{0\}, \ \text{and} \ r >0.
\end{equation*}
We easily see that 
\begin{equation} \label{eq-calpha}
  \left| K_\alpha(x) \right| \le \frac{c_\alpha}{|x|^{\alpha}}, 
  \quad \text{for every} \  x \in \mathbb{R}^n \setminus \{0\},
\end{equation}
where $c_\alpha := \max_{|x| = 1} \left| K_\alpha(x) \right|$. 
%\begin{remark} \label{rem-p-alpha-n}
In particular, $K_\alpha$ belongs to the space $L^q_{\rm weak}(\R^n)$, for $q = \dfrac{n}{\alpha}$, that is
\begin{equation} \label{eq-Lpw}
  \text{\rm meas} \big\{ x \in \R^n \mid | K_\alpha (x) | \ge \lambda  \big\} \le \left( \frac{C}{\lambda}\right)^q, 
\quad \text{for every} \quad \lambda >0,
\end{equation}
for some non-negative constant $C$. Here {\rm meas} $E$ denotes the Lebesgue measure of the set $E$. Moreover we 
define the \emph{seminorm} of $K_\alpha$ as follows
\begin{equation*} %\label{eq-Lpw-sn}
  \left\| K_\alpha \right\|_{L^q_{\rm weak}(\R^n)} := \inf \big\{ C \ge 0 \mid \text{\eqref {eq-Lpw} holds} \big\}.
\end{equation*}
From \eqref{eq-calpha} it plainly follows that $C \le c_\alpha \omega_n^{\alpha/n}$. We next recall two elementary 
inequalities that will be useful in the sequel. For every $R>0$ we have that:
  \begin{itemize}
    \item $K_\alpha \in L^q\big(\big\{ x \in \R^n \mid |x| \le R \big\}\big)$ if, and only if, $q < \dfrac{n}{\alpha}$. 
Moreover, there exists a positive constant $c_{\alpha,q}$, only depending on $K_\alpha, n$ and $q$, such that 
\begin{equation}\label{eq-pa>n}
  \left\| K_\alpha \right\|_{L^q(\{ x \in \R^n \mid |x| \le R\})} \le 
  c_{\alpha,q} R^{\frac{n}{q} - \alpha},
\end{equation}
    \item $K_\alpha \in L^q\big(\big\{ x \in \R^n \mid |x| \ge R \big\}\big)$ if, and only if, $q > \dfrac{n}{\alpha}$. 
Moreover, there exists a positive constant $c_{\alpha,q}$, only depending on $K_\alpha, n$ and $q$, such that 
\begin{equation}\label{eq-pa<n}
  \left\| K_\alpha \right\|_{L^q(\{ x \in \R^n \mid |x| \ge R \})} \le 
  c_{\alpha,q} R^{\frac{n}{q} - \alpha}.
\end{equation}
  \end{itemize}
%\end{remark}
The following weak Young inequality holds (see Theorem 1, p. 119 in \cite{undici}, where this result is referred to 
as Hardy-Littlewood-Sobolev theorem for fractional integration).
\begin{theorem} \label{Th-weak-Young}
Let $K_\alpha$ be a continuous  homogeneous function of degree $- \alpha$, with $0 < \alpha < n$. Let $p, q$ be such 
that $1 \le p < q < + \infty$ and that $1 + \dfrac{1}{q} = \dfrac{1}{p} + \dfrac{\alpha}{n}$. Then, for every $f \in 
L^p(\R^n)$ the integral $K_\alpha * f (x)$ is convergent for almost every $x \in \R^n$. Moreover, 

\noindent - if $p>1$, then there exists $C_{\alpha, p} >0$ such that $\norm{K_\alpha * f}_{L^q (\mathbb{R}^n)} \leq 
C_{\alpha, p} \norm{f}_{L^p (\mathbb{R}^n)}$,
%  \begin{equation*}
% \norm{K_\alpha * f}_{L^q (\mathbb{R}^n)} \leq C_{\alpha, p} \norm{f}_{L^p (\mathbb{R}^n)}, 
% \quad \text{for every} \quad f\in L^p (\mathbb{R}^n);
% \end{equation*}

\noindent - if $p=1$, then there exists $C_{\alpha,1} >0$ such that $\left\| K_\alpha * f \right\|_{L^q_{\rm 
weak}(\R^n)} \le C_{\alpha, 1} \norm{f}_{L^1 (\mathbb{R}^n)}$.
% $K_\alpha * f \in L^q_{\rm weak}(\R^n)$, and
% there exists $C_{\alpha, 1} >0$ such that  
% \begin{equation*}
%   \text{\rm meas} \big\{ x \in \R^n \mid | K_\alpha * f | \ge \lambda  \big\} \le \left( \frac{C_{\alpha, 1} 
% \norm{f}_1}{\lambda}\right)^q, \quad \text{for every} \quad \lambda >0. 
% \end{equation*}
% Here {\rm meas} $E$ denotes the Lebesgue measure of the set $E$.
\end{theorem}

In order to prove the Morrey embedding \eqref{eq-Morrey} and the compactness estimate \eqref{eq-Sobolev-unif}, we 
state and prove the following lemma. 

\begin{lemma} \label{14}
Let $K_\alpha \in C^1(\mathbb{R}^n\backslash \lbrace0\rbrace)$ be any homogeneous function of degree $- \alpha$, with 
$0<\alpha<n$. Then there exists a positive constant $M_\alpha$ such that 
\[
\abs{K_\alpha(x)-K_\alpha(y)}\leq M_\alpha \frac{\abs{x-y}}{\abs{x}^{\alpha+1}}, \quad \text{for every} \quad 
x,y\in\mathbb{R}^n\backslash \lbrace 0 \rbrace \ \text{such that} \ \abs{x-y}\leq \dfrac{\abs{x}}{2}.
\]
\end{lemma}
\noindent \textit{Proof.} We first prove the result for $x$ such that $\abs{x}=1$. In this case 
$\abs{x-y}\leq\dfrac{1}{2}$ 
and by the Mean Value Theorem there exists $\theta\in(0,1)$ such that
\begin{equation*} %\label{MVT}
|K_\alpha(x)-K_\alpha(y)|=|\scalare{(x-y)}{\nabla K_\alpha(\theta x+(1-\theta y))}|\leq M_\alpha\abs{x-y},
\end{equation*}
where
\begin{equation} \label{Emme}
M_\alpha=\max_{\frac{1}{2}\leq |z| \leq \frac{3}{2}}{\abs{\nabla K_\alpha(z)}}.
\end{equation}
% due to the fact that
% \[
% \dfrac{1}{2}\leq \abs{\theta x+(1-\theta)y}\leq\dfrac{3}{2}
% \]
Consider now a general choice of $x,y\in\mathbb{R}^n\backslash \lbrace0\rbrace$ with $\abs{x-y}\leq\dfrac{\abs{x}}{2}$. 
Being $K_\alpha$ homogeneous of degree $- \alpha$, we obtain
\[
\abs{K_\alpha(x)-K_\alpha(y)}=\frac{1}{\abs{x}^{\alpha}} \left 
|K_\alpha\left(\frac{x}{|x|}\right)-K_\alpha\left(\frac{y}{|x|}\right) \right | 
% \stackrel{\eqref{MVT}}
{\leq} \frac{M_\alpha}{\abs{x}^{ \alpha}}\left |\frac{x}{|x|} - \frac{y}{|x|} \right |= 
M_\alpha\frac{\abs{x-y}}{\abs{x}^{\alpha+1}},
\]
being $M_\alpha$ as in \eqref{Emme}, because $\left| \dfrac{x}{|x|}\right|= 1$. $\Box$

\medskip

In order to prove the Morrey embedding \eqref{eq-Morrey} and the compactness of the Sobolev embedding 
for $p < q < p^*$ we rely on the following argument. We choose any $u \in C_0^\infty(\R^n), h \in \R^n$ and 
we set
\begin{equation} \label{eq-vx}
  v(x) := u(x + h) - u(x), \qquad \text{for every} \quad x \in \R^n, 
\end{equation}
then
\begin{eqnarray}
%\abs{u(x+h)-u(x)} 
v(x) &=&\integrale{\lbrace y\in\mathbb{R}^n:\abs{x+h-y} \geq 
2\abs{h}\rbrace}{} {\langle \nabla\Gamma(x+h-y)-\nabla\Gamma(x-y), \nabla u(y) \rangle}{y} \nonumber \\ 
& + & \integrale{\lbrace y\in\mathbb{R}^n:\abs{x+h-y} < 2\abs{h}\rbrace}{}{\langle{\nabla\Gamma(x+h-y)},{\nabla 
u(y)} \rangle} {y} \nonumber\\
& + & \integrale{\lbrace y\in\mathbb{R}^n:\abs{x+h-y} < 2\abs{h}\rbrace}{}{- \langle {\nabla\Gamma(x-y)}, {\nabla 
u(y)}\rangle}{y} =: I_A (x) + I_B (x)  + I_C (x). \label{usplit} 
\end{eqnarray}
We next rely on Lemma \ref{14} and on \eqref{eq-grad-fund-sol-lapl} to estimate the terms $I_A, I_B$ and $I_C$ as 
follows
\begin{equation} \label{eq-bound-ABC}
\begin{split}
 | I_A (x) | & \le M |h|   \integrale{\lbrace y\in\mathbb{R}^n:\abs{x+h-y} \geq 
2\abs{h}\rbrace}{}{\frac{1}{\abs{x+h-y}^n}\abs{\nabla u(y)}}{y},  \\ %\label{eq-bound-A}
 | I_B (x) | & \le c_n  \integrale{\lbrace y\in\mathbb{R}^n:\abs{x+h-y} < 
2\abs{h}\rbrace}{}{\frac{1}{\abs{x+h-y}^{n-1}}\abs{\nabla u(y)}}{y}, \\ %\label{eq-bound-B} 
 | I_C (x) | & \le c_n  \integrale{\lbrace y\in\mathbb{R}^n:\abs{x+h-y} < 
2\abs{h}\rbrace}{}{\frac{1}{\abs{x-y}^{n-1}}\abs{\nabla u(y)}}{y},  %\label{eq-bound-C}
\end{split}
\end{equation}
where $M := \max\limits_{\frac{1}{2}\leq |z| \leq \frac{3}{2}, j,k = 1 , \dots, n}{\abs{ \partial^2_{x_j x_k} 
\Gamma(z)}}$.

\subsection{The Sobolev and Morrey embedding theorems} 

As we said in the Introduction, Theorem \ref{Th-weak-Young} combined with \eqref{eq-conv-grad} immediately yields the 
following result
\begin{theorem} \label{Th-Sobolev}
Let $1<p<n$. Then there exists $C_p>0$ such that:
\begin{equation*}
\norm{u}_{p^*}\leq C_p \norm{\nabla u}_p  \quad \text{for every} \quad  u \in L^1_{\text {\rm 
loc}}(\mathbb{R}^n) \quad  \text{such that} \quad \nabla u \in L^{p} (\mathbb{R}^n),
\end{equation*}
where $p^*=\dfrac{np}{n-p}$.
\end{theorem}

%\subsection{The Morrey embedding} 

% then, by \eqref{usplit}, we write 
% \begin{equation} \label{eq-vsplit}
%   v(x) = I_A(x) + I_B(x) + I_C(x), \qquad \text{for every} \quad x \in \R^n. 
% \end{equation}
% 
% We fix $u \in C_0^\infty(\R^n), h \in \R^n$ and we set
% \begin{equation} \label{eq-vx}
%   v(x) := u(x + h) - u(x), \qquad \text{for every} \quad x \in \R^n, 
% \end{equation}

We next turn our attention on the Morrey's Theorem.

\begin{theorem}[Morrey's Theorem] \label{e}
Let $u:\mathbb{R}^n\rightarrow\mathbb{R}$ and $p>n$. If $\nabla u \in \sob{p}$, then $u$ is continuous and 
\begin{equation} \label{Morrey-ineq}
\abs{u(x+h)-u(x)}\leq  C_{n,p} \norm{\nabla u}_p\abs{h}^{1-\frac{n}{p}}, \qquad \text{for every} \quad x, h \in \R^n,
\end{equation}
for some positive constant $C_{n,p}$ depending only on $p$ and $n$.

In particular $\sobolev{p}$ is continuously embedded in the space of H\"older continuous 
functions $C^\beta(\mathbb{R}^n)$, with $\beta=1-\frac{n}{p}$. 
\end{theorem}

% $u\in C^\alpha(\mathbb{R}^n)$, with
% \[
% \alpha=1-\frac{n}{p}
% \]
% and fixed $x\in\mathbb{R}^n$, for all $h\in\mathbb{R}^n$ it holds
% \begin{equation}
% \label{Morrey-ineq}
% \abs{u(x+h)-u(x)}\leq \bmicr C_{n,p} \emicr \norm{\nabla u}_p\abs{h}^\alpha
% \end{equation}
% with a positive constant $C_{n,p}>0$ depending on $p$ and $n$.

\noindent \textit{Proof.} Consider a function $u \in C_0^\infty(\mathbb{R}^n)$, let $x, h \in \R^n$, and $v(x)$ be 
the function defined in \eqref{usplit}. We next estimate  $I_A, I_B$ and $I_C$ by using the H\"older inequality.

From the first inequality in \eqref{eq-bound-ABC} and \eqref{eq-pa<n}, with $\alpha = n$ and $q = \dfrac{p}{p-1}$, we 
obtain
\begin{equation*} %\label{eq-Mor-A}
\begin{split}
    | I_A (x) | \le & \, M |h| \norm{\nabla u}_{L^p(\mathbb{R}^n)} 
    \left\|{\frac{1}{|z|^n}}\right\|_{L^q(\{ |z| \ge 2 |h|\})} \\
    = & \, M |h| \norm{\nabla u}_{L^p(\mathbb{R}^n)} c_{n,q} (2|h|)^{\frac{n(p-1)}{p} - n} = 
    M_A\norm{\nabla u}_{L^p(\mathbb{R}^n)}\abs{h}^{1-\frac{n}{p}},
\end{split}
\end{equation*}
for some positive constant $M_A$ only depending on $M, p$ and $n$. Moreover, from the second inequality in 
\eqref{eq-bound-ABC} and \eqref{eq-pa>n}, with $\alpha = n-1$ and $q = \dfrac{p}{p-1}$, 
we find
\begin{equation*} %\label{eq-Mor-B}
\begin{split}
    | I_B (x) | \le & \, c_n \norm{\nabla u}_{L^p(\mathbb{R}^n)} 
    \left\|{\frac{1}{|z|^{n-1}}}\right\|_{L^q( \{|z| \le 2 |h|\})} \\
    = & \, c_n \norm{\nabla u}_{L^p(\mathbb{R}^n)} c_{n-1,q} (2|h|)^{\frac{n(p-1)}{p} - n+1} = 
    M_B \norm{\nabla u}_{L^p(\mathbb{R}^n)}\abs{h}^{1-\frac{n}{p}},
\end{split}
\end{equation*}
where $M_B$ is a positive constant only depending on $p$ and $n$. Finally the last term in \eqref{usplit} can be 
estimated similarly to the second one, observing that $\abs{x-y}\leq \abs{x+h-y}+\abs{h}< 3\abs{h}$, thus getting
\begin{equation*} %\label{eq-Mor-B}
\begin{split}
    | I_C (x) | \le & c_n \norm{\nabla u}_{L^p(\mathbb{R}^n)} 
    \left\|{\frac{1}{|z|^{n-1}}}\right\|_{L^q( \{|z| \le 3 |h|\})} = M_C 
\norm{\nabla u}_{L^p(\mathbb{R}^n)}\abs{h}^{1-\frac{n}{p}},
\end{split}
\end{equation*}
for some positive constant $M_C$. We note that  the $L^q$ norms of the functions $z \mapsto |z|^{-n}$ and $z 
\mapsto |z|^{-n+1}$ appearing in the above estimates are finite thanks to the assumption $p > n$.
Then \eqref{Morrey-ineq} is obtained with $C_{n,p} := M_A + M_B + M_C$.  This proves our claim for $u \in 
C_0^\infty(\R^n)$. The general case follows by a density argument. \hfill $\Box$

\medskip

We next prove the compactness of the Sobolev embedding \eqref{eq-Sobolev-unif} for $p < q < p^*$ starting again from
\eqref{usplit}. Here we use the Young inequality instead of the H\"older inequality.
%More specifically, the following result holds true.  

\begin{theorem} \label{Uniform continuity}
Let $p, q \ge 1$ be such that $p < q < p^* := \frac{np}{n-p}$. Then there exists a positive constant $C_{p,q}$ 
depending on $n, p, q$ such that 
\[
\norm{u(h + \cdot)-u}_{L^q(\mathbb{R}^n)}\leq  C_{p,q} \, \|\nabla u\|_{L^p(\mathbb{R}^n)} 
\, \abs{h}^{n \left(\frac{1}{q} - \frac{1}{p^*}\right)}. 
\]
for every $u \in C^{\infty}_0(\mathbb{R}^n)$ and for any $h \in \mathbb{R}^n$. In particular, as long as $q < p^*$, we 
have 
\[
\norm{u(h + \cdot)-u}_{L^q(\mathbb{R}^n)} \rightarrow 0 \qquad \textnormal{as} \qquad |h| \rightarrow 0.
\]
\end{theorem}

\noindent \textit{Proof.} Consider a function $u \in C_0^\infty(\mathbb{R}^n)$, let $x, h \in \R^n$, and $v(x)$ be 
the function defined in \eqref{usplit}. We next estimate the $L^q$ norm of $I_A, I_B$ and $I_C$ by using the Young
inequality. To this aim we introduce the exponent $r$ defined by the identity
\begin{equation} \label{esponente-r}
1 + \frac{1}{q} = \frac{1}{r} + \frac{1}{p},
%r := \frac{1}{\frac{1}{q} + 1 - \frac{1}{p}};
\end{equation}
and we note that 
\begin{equation} \label{condizionesur}
1 < r < \frac{n}{n-1} \,\, \Leftrightarrow \,\, p < q < p^*.
\end{equation}
From the first inequality in \eqref{eq-bound-ABC} and \eqref{eq-pa<n}, with $\alpha = n$, we obtain
\begin{equation*} %\label{eq-Mor-A}
\begin{split}
    \left\| I_A \right\|_{L^q(\R^n)} \le & \, M |h| \norm{\nabla u}_{L^p(\mathbb{R}^n)} 
    \left\|{\frac{1}{|z|^n}}\right\|_{L^r(\{ |z| \ge 2 |h| \})} \\
    = & \, M |h| \norm{\nabla u}_{L^p(\mathbb{R}^n)} c_{n,r} (2|h|)^{\frac{n}{r} - n} 
    = C_A(n,r) |h|^{n\left(\frac{1}{q}-\frac{1}{p^*}\right)},
\end{split}
\end{equation*}
From the second inequality in \eqref{eq-bound-ABC} and \eqref{eq-pa>n}, with $\alpha = n-1$, we obtain
\begin{equation*} %\label{eq-Mor-B}
\begin{split}
    \left\| I_B \right\|_{L^q(\R^n)} \le & \, c_n \norm{\nabla u}_{L^p(\mathbb{R}^n)} 
    \left\|{\frac{1}{|z|^{n-1}}}\right\|_{L^r(\{ |z| \le 2 |h| \})} \\
    = & \, c_n \norm{\nabla u}_{L^p(\mathbb{R}^n)} c_{n-1,r} (2|h|)^{\frac{n}{r} - n+1} = 
 C_B(n,r) \norm{\nabla u}_{L^p(\mathbb{R}^n)}\abs{h}^{n\left(\frac{1}{q}-\frac{1}{p^*}\right)},
\end{split}
\end{equation*}
where $C_B(n,r)$ is a constant depending on $n$ and $r$ (and thus on $n, p, q$). The same argument applies to $I_C$, so 
that, provided that we consider the norm $\left\|{\frac{1}{|z|^{n-1}}}\right\|_{L^r( \{|z| \le 3 |h|\})}$ instead of 
$\left\|{\frac{1}{|z|^{n-1}}}\right\|_{L^r(\{ |z| \le 2 |h|\})}$, as in the proof of the Morrey's theorem. We then find
\begin{equation*} %\label{eq-Mor-B}
\begin{split}
    \left\| I_C \right\|_{L^q(\R^n)} \le C_C(n,r) \norm{\nabla 
     u}_{L^p(\mathbb{R}^n)}\abs{h}^{n\left(\frac{1}{q}-\frac{1}{p^*}\right)},
\end{split}
\end{equation*} 
We note that  the $L^r$ norms of the functions $z \mapsto |z|^{-n}$ and $z \mapsto |z|^{-n+1}$ appearing in the 
above estimates are finite if, and only if, the condition \eqref{condizionesur} is satisfied. The thesis is 
obtained with $C_{p,q} := C_A(n,r) + C_B(n,r) + C_C(n,r)$ for $u \in C_0^\infty(\R^n)$. The 
general case follows by a density argument. \hfill $\Box$

\subsection{A more general compactness result}

We note that the Theorem \ref{Uniform continuity} only applies to a kernel that is homogeneous of degree $-n+1$. 
Actually, the method used in its proof also applies to any general homogeneous kernel $K_\alpha$, with $0 < \alpha < 
n$, as those considered in Theorem \ref{Th-weak-Young}. In the following statement we denote by $u$ the convolution 
$K_{\alpha} * f$, that is
\[
u(x)=\integ{K_{\alpha} (x-y) f(y)}{y}.
\]

\begin{theorem} \label{Uniform continuity-bis}
Let $K_\alpha$ be a $C^1(\R^n \setminus \{ 0 \})$ homogeneous function of degree $- \alpha$, with $1 < \alpha < n$, and 
let $p, q \ge 1$ be such that $q > p$ and 
% p < q < + \infty$ and $1 + \dfrac{1}{q} > \dfrac{1}{p} + \dfrac{\alpha}{n}
\begin{equation} \label{eq-pqalpha}
  1 - \dfrac{\alpha+1}{n} < \dfrac{1}{p} - \dfrac{1}{q} <  1 - \dfrac{\alpha}{n}.
\end{equation}
Then there exists a positive constant $\widetilde{C}_{p,q}$, depending on $n, p, q$, such that 
\[
\norm{u(h + \cdot)-u}_{L^q(\mathbb{R}^n)}\leq  \widetilde{C}_{p,q} \, \| f \|_{L^p(\mathbb{R}^n)} 
\, \abs{h}^{\frac{n}{r} - \alpha},
\]
for every $f \in \sobb{p}{{\mathbb{R}^{n}}}$ and $h \in \R^n$. Here $r$ is the constant introduced in 
\eqref{esponente-r}, that is $1 + \dfrac{1}{q} = \dfrac{1}{r} + \dfrac{1}{p}$. Moreover the exponent 
$\dfrac{n}{r} - \alpha$ is strictly positive.
\end{theorem}

\noindent \textit{Proof.} We choose $x, h \in \R^n$ and let $v$ be defined as in \eqref{eq-vx}: $v(x) = u(x+h) - u(x)$, 
and we consider three integrals $v(x) = {I}_A(x) + {I}_B(x) + {I}_C(x)$ as in \eqref{usplit}. 
% Then
% \begin{eqnarray}
% \abs{u(x+h)-u(x)} &=&\integrale{\lbrace y\in\mathbb{R}^n:\abs{x+h-y} \geq 
% 2\abs{h}\rbrace}{}{\abs{K_{\alpha}(x+h-y)-K_{\alpha}(x-y)}\abs{f(y)}}{y} \nonumber \\ 
% & + & \integrale{\lbrace y\in\mathbb{R}^n:\abs{x+h-y} < 2\abs{h}\rbrace}{}{\abs{K_{\alpha}(x+h-y)}\abs{f(y)}}{y} 
% \nonumber\\
% & + & \integrale{\lbrace y\in\mathbb{R}^n:\abs{x+h-y} < 2\abs{h}\rbrace}{}{\abs{K_{\alpha}(x-y)}\abs{f(y)}}{y} =: 
% \tilde{I}_A + \tilde{I}_B + \tilde{I}_C \label{usplittilde}. 
% \end{eqnarray}
% Then,  
% \begin{equation} \label{eq-vsplit-tilde}
%   v(x) = \tilde{I}_A(x) + \tilde{I}_B(x) + \tilde{I}_C(x), \qquad \text{for every} \quad x \in \R^n
% \end{equation}
We proceed as we did in the proof of Theorem \ref{e}. We find 
\begin{equation*}
\begin{split}
 | I_A (x) | & \le M_\alpha |h|   \integrale{\lbrace y\in\mathbb{R}^n:\abs{x+h-y} \geq 
2\abs{h}\rbrace}{}{\frac{1}{\abs{x+h-y}^{\alpha + 1}}\abs{f(y)}}{y}, \\
 | I_B (x) | & \le c_\alpha  \integrale{\lbrace y\in\mathbb{R}^n:\abs{x+h-y} < 
2\abs{h}\rbrace}{}{\frac{1}{\abs{x+h-y}^{\alpha}}\abs{f(y)}}{y}, \\
 | I_C (x) | & \le c_\alpha  \integrale{\lbrace y\in\mathbb{R}^n:\abs{x+h-y} < 
2\abs{h}\rbrace}{}{\frac{1}{\abs{x-y}^{\alpha}}\abs{f(y)}}{y}, 
\end{split}
\end{equation*}
being $M_\alpha := \max\limits_{\frac{1}{2}\leq |z| \leq \frac{3}{2}}{\abs{ \nabla K_\alpha(z)}}$. In order to use the 
Young inequality, we recall that 
\begin{equation*}
%\begin{split}
\left\| \frac{1}{\abs{z}^{\alpha + 1}} \right\|_{L^r(\lbrace y\in\mathbb{R}^n:\abs{z} \geq 2\abs{h}\rbrace)}  
%& \bigg( \integrale{\{\abs{z}\geq 2\abs{h}\}}{}{\frac{1}{\abs{z}^{r(\alpha + 1)}}}{x} 
%\bigg)^{\frac{1}{r}}\\
%=& \bigg( \integrale{2\abs{h}}{+\infty}{\omega_n\rho^{n-1+r(\alpha+1)}}{\rho} \bigg)^{\frac{1}{r}} 
= c_{\alpha+1, r} (2|h|)^{\frac{n}{r} - (\alpha + 1)},
%\end{split}
\end{equation*}
where $r$ is the exponent introduced in \eqref{esponente-r} and 
%$\widetilde{C}_A(n,r)$ is a constant depending on $n$ and $r$  (and so on $n, p, q$). Note 
we note that the above integral converges if, and only if, $({\alpha + 1})r > {n}$. Moreover, 
\begin{equation*}
%\begin{split}
\left\| \frac{1}{\abs{z}^{\alpha}} \right\|_{L^r(\lbrace z \in\mathbb{R}^n:\abs{z} \leq 2\abs{h}\rbrace)}  
%& \bigg( \integrale{\{\abs{z}\geq 2\abs{h}\}}{}{\frac{1}{\abs{z}^{r(\alpha + 1)}}}{x} 
%\bigg)^{\frac{1}{r}}\\
%=& \bigg( \integrale{2\abs{h}}{+\infty}{\omega_n\rho^{n-1+r(\alpha+1)}}{\rho} \bigg)^{\frac{1}{r}} 
= c_{\alpha, r} (2|h|)^{\frac{n}{r} - \alpha},
%\end{split}
\end{equation*}
and the above integral converges if, and only if, ${\alpha} r <{n}$. Summarizing, the two above integrals are finite 
if, and only if,
\begin{equation}
\label{cond-su-r}
\frac{n}{\alpha + 1} < r < \frac{n}{\alpha}, 
\end{equation}
which is equivalent to \eqref{eq-pqalpha}. We also note that the exponent $\dfrac{n}{r} - \alpha$ is strictly positive.

We next proceed as in the proof of Theorem \ref{e}. By the Young inequality we deduce
\begin{eqnarray*}
& \left\| I_A \right\|_{L^q(\R^n)} \le \widetilde{C}_A(n,\alpha, r) |h|^{\left(\frac{n}{r}- \alpha\right)} 
\left\| f \right\|_{L^p(\R^n)}, \\
& \left\| I_B \right\|_{L^q(\R^n)} \le \widetilde{C}_B(n,\alpha, r) |h|^{\left(\frac{n}{r}- \alpha\right)} 
\left\| f \right\|_{L^p(\R^n)}, \\
& \left\| I_C \right\|_{L^q(\R^n)} \le \widetilde{C}_C(n,\alpha, r) |h|^{\left(\frac{n}{r}- \alpha\right)} 
\left\| f \right\|_{L^p(\R^n)}, 
\end{eqnarray*}
where $\widetilde{C}_A(n,\alpha, r),  \widetilde{C}_B(n,\alpha, r)$, and $\widetilde{C}_C(n,\alpha, r)$ are positive 
constants only depending on $n, p, q$, on $c_\alpha$ and on $M_\alpha$. The thesis is obtained with 
$\widetilde{C}_{p,q} 
:= \widetilde{C}_A(n,\alpha, r) +  \widetilde{C}_B(n,\alpha, r) + \widetilde{C}_C(n,\alpha, r)$. \hfill $\Box$

\begin{remark} \label{listona}
The statement of Theorem \ref{Uniform continuity-bis} is more involved with respect to that of Theorem 
\ref{Uniform continuity} as we don't have a natural counterpart of the Sobolev exponent $p^*$ for any $\alpha \in 
]0,n[$. We list here the explicit conditions on $q$ for the validity of \eqref{eq-pqalpha}. We discuss 
specifically the case $1 < p < n$ as we are interested 
in the compactness of the Sobolev embedding. 
  \begin{itemize}
\item[(i)] If \fbox{$0 < \alpha < n-2$} then
\begin{itemize}
\item[(i.1)] If $\displaystyle 1 < p < \frac{n}{n - \alpha}$ then 
\[  
\frac{np}{n - p (n - 1 - \alpha)} < q < \frac{np}{n - p (n - \alpha)}
\]
\item[(i.2)] If $\displaystyle \frac{n}{n - \alpha} < p < \frac{n}{n - 1 - \alpha}$ then 
\[ 
q > \max \left \{1, \frac{np}{n - p (n - 1 - \alpha)} \right \} = \frac{np}{n - p (n - \alpha)} 
\]
\item[(i.3)] If $\displaystyle  \frac{n}{n - 1 - \alpha} < p < n$ then no values of $q$ are available
\end{itemize}
\item[(ii)] If \fbox{$\alpha = n - 2$} then
\begin{itemize}
\item[(ii.1)] If $\displaystyle 1 < p < \frac{n}{2}$ then 
\[  
p^* = \frac{np}{n - p} < q < \frac{np}{n - 2p}
\]
\item[(ii.2)] If $\displaystyle \frac{n}{2} < p < n$ then 
\[ 
q > \max \left \{1, \frac{np}{n - p} \right \} = \frac{np}{n - p} = p^*
\]
\end{itemize}
\item[(iii)] If \fbox{$n-2 < \alpha < n-1$} then 
\begin{itemize}
\item[(iii.1)] If $\displaystyle 1 < p < \frac{n}{n - \alpha}$ then 
\[  
\frac{np}{n - p (n - 1 - \alpha)} < q < \frac{np}{n - p (n - \alpha)}
\]
\item[(iii.2)] If $\displaystyle \frac{n}{n - \alpha} < p < n$ then 
\[ 
q > \max \left \{1, \frac{np}{n - p (n - \alpha)} \right \} = \frac{np}{n - p (n - 1 - \alpha)} 
\]
\end{itemize}
\item[(iv)] If \fbox{$\alpha = n - 1$} then for any $1 < p < n$
\[  
p < q < p^* = \frac{np}{n - p}
\]
\item[(v)] If \fbox{$n - 1 < \alpha < n$} then for any $1 < p < n$
\[  
\frac{np}{n - p (n - 1 - \alpha)} < q < \frac{np}{n - p (n - \alpha)}.
\]
\end{itemize}
\end{remark}

\section{Continuous and compact embeddings for degenerate Kolmogorov equations}
\setcounter{section}{3} \setcounter{equation}{0}
\setcounter{theorem}{0} 

The aim of this section is to prove compactness estimates for weak solutions to a family of degenerate Kolmogorov
operators that includes the one in \eqref{e-Kolmo-1} as the simplest prototype. Specifically, we consider 
operators in this form
\begin{equation} \label{ultraparabolica}
\L u(x,t) :=\sum_{i,j=1}^{m_0}{\partial_{x_i}(a_{i,j}(x,t)\partial_{x_j}u(x,t))}+
\sum_{i,j=1}^{N}{b_{i,j}x_i\partial_{x_j} u(x,t)}-\partial_t u(x,t),
\end{equation}
with $(x,t) =(x_1,...,x_N,t) \in \rnn$. Here $m_0\in\mathbb{N}$ is such that $1\leq m_0\leq N$. In the 
sequel we will also use the following notation $z:=(x,t)$, and we always assume the following hypotheses:
\begin{itemize}
\item[(H.1)] The matrix $A=(a_{i,j}(z))_{i,j=1,...,m_0}$ is symmetric, with measurable coefficients and there exists a 
positive constant  $\mu$ such that
\[
\mu^{-1}\abs{\xi}^2 \leq \sum_{i,j=1}^{m_0}{a_{i,j}(z)\xi_i\xi_j} \leq \mu\abs{\xi}^2
\]
for all $z\in\rnn$ and $\xi\in\mathbb{R}^{m_0}$.
\item[(H.2)] The matrix $B$ has constant coefficients. Moreover there exists a basis of $\rn$ such that the 
matrix $B$ can be written in a canonical form:
\begin{equation*}
B=\begin{pmatrix}
0 & B_1 & 0 & ... & 0 \\
0 & 0 & B_2 &... & 0 \\
\vdots & \vdots &  \vdots & \ddots & \vdots \\
0 & 0 & 0 &... & B_r \\
0 & 0 & 0 &... & 0
\end{pmatrix}
\end{equation*}
where $B_k$ is a matrix $m_{k-1}\times m_k$ with rank $m_k$, $k=1,2,...,r$ with
\[
m_0\geq m_1\geq ... \geq m_r\geq 1, \ \ \ \mathrm{e} \ \ \ \sum_{k=0}^r m_k=N.
\]
\end{itemize}
We prove in detail Theorem \ref{main-th-1} along the same techniques outlined in the Introduction for the operator $\L$ 
in \eqref{e-Kolmo-1}. In particular, we will rely on a representation formula analogous to \eqref{eq-fund-sol2} in 
terms 
of the fundamental solution to the operator
\begin{equation} \label{ultraparabolica-0} 
  \L_0 u :=\sum_{i=1}^{m_0}{\partial_{x_i}^2 u}+\sum_{i,j=1}^{N}{b_{i,j}x_i\partial_{x_j}u}-\partial_t u.
\end{equation}
\textit{Remark 1}.
It is know that Assumption (H.2) is equivalent to the assumption of hypoellipticity of $\L_0$ (see \cite{dodici} and 
its bibliography). This means that any function  $u$ which is a distributional solution to $\L_0 u=f$ in some open 
subset $\Omega$ of $\rnn$ is a $C^\infty$ function whenever $f$ is $C^\infty$.

In order to simplify our statements, in the sequel we adopt the following compact notation. If $I_{m_0}$ is the 
identity matrix 
$m_0\times m_0$, we set
\[
\Delta_{m_0}=\sum_{i=1}^{m_0}\partial_{x_i}^2, \quad \quad 
Y=\sum_{i,j=1}^{N}{b_{i,j}x_i\partial_{x_j}}-\partial_t,
\quad \quad 
A_0=\begin{pmatrix}
I_{m_0} & 0 \\
0 & 0 \\
\end{pmatrix}.
\]
In particular we will write
\[
\L_0 =\Delta_{m_0}+Y =\mathrm{div}(A_0\nabla)+Y.
\]
We next introduce the non-Euclidean geometric setting suitable for the study of $\L$, the fundamental solution of 
$\L_0$ and the definition of the convolution with homogeneous kernels. 

\subsection{Dilation and translation groups associated to $\L$}
% Let us consider the second order partial differential equation \eqref{ultraparabolica}
% \begin{equation*}
% 

% \end{equation*}
% where:
% \begin{itemize}
% \item $z=(x,t)=(x_1,...,x_n,t) \in \rnn$.
% \item $m_0\in\mathbb{N}$ is such that $1\leq m_0\leq n$.
% \end{itemize}
% With the aim to deal with equation (\ref{ultraparabolica}) we will need some results concerning the operator $L_0$ 
% defined as follows:
We recall here some invariance properties of the operator $\L_0$. We refer to \cite{dodici} where the definition of 
translation group and dilation group for Kolmogorov operators have been given for the first time. Let
\begin{equation} \label{eq-translation}
(x,t)\circ (\xi,\tau)=(\xi+E(\tau)x,t+\tau), \ \ \ E(t)=\exp(-tB^T), \ (x,t),(\xi,\tau)\in\rnn
\end{equation}
\begin{equation} \label{eq-dilation}
  D(\lambda)=\mathrm{diag}(\lambda I_{m_0},\lambda^3 I_{m_1},...,\lambda^{2r+1} I_{m_r},\lambda^2), \ \ \ \lambda>0
\end{equation}
where $I_{m_j}$ denotes the identity matrix $m_j\times m_j$. It is known that $(\rnn,\circ)$ is a 
\emph{non commutative} group, and $\L_0$ is invariant with respect to the \emph{left translations} of 
$(\rnn,\circ)$, in the following sense: if we choose any $\zeta \in \rnn$ and we set $v(z) := u (\zeta \circ z)$ and  
$g(z) := f (\zeta \circ z)$, then the have
\[
\L_0 u (z) = f (z) \quad \Leftrightarrow \quad \L_0 v (z) =g (z).
\]
Moreover, $\L_0$ is invariant with respect to $(D(\lambda))_{\lambda>0}$, with the following meaning:
\[
\L_0 u (z) = f (z) \quad \Leftrightarrow \quad \L_0 w (z) = \lambda^2 h (z).
\]
Now $w(z) := u(D(\lambda) z)$, $h(z) := f(D(\lambda) z)$, and $\lambda$ is any positive constant. The zero of the group 
is $(0,0)$ and the inverse of $(\xi, \tau)$ is $(\xi, \tau)^{-1} = (- E(\tau)\xi, - \tau)$. Moreover the following 
\emph{distributive property} holds:
\[
   D(\lambda) \left( z \circ \zeta \right) = \left( D(\lambda) z \right) \circ \left( D(\lambda) \zeta \right) 
\]
We summarize the above properties by saying that $\L_0$ is invariant with respect to the \emph{homogeneous group} 
$\left(\rnn,\circ, (D(\lambda))_{\lambda>0} \right)$. 

In the sequel we will use the following notation 
\[
   D(\lambda) = \text{diag} (\lambda^{\alpha_1}, \lambda^{\alpha_2}, \dots , \lambda^{\alpha_N}, \lambda^{2}),
\]
where, in accordance with \eqref{eq-dilation} $\alpha_1 = \alpha_2 = \dots = \alpha_{m_0} = 1$,  $\alpha_{m_0 + 1} = 
\dots = \alpha_{m_0 + m_1} = 3, \dots , \alpha_{m_0 + m_1 + \dots + m_{r-1} + 1} = \dots = \alpha_{N} = 2 r +1$.
%\subsection{The homogeneous norm with respect to the dilation group} 

We define now a norm of $\rnn$ homogeneous of degree one with respect to the dilation introduced before.
\begin{definition} \label{def-norm}
\label{normaomogenea}
For all $z\in\rnn\setminus \{0\}$, we define the norm $\norm{z}=\rho$, as the unique positive solution to:
\[
\frac{x_1^2}{\rho^{2\alpha_1}}+\frac{x_2^2}{\rho^{2\alpha_2}}+...+\frac{x_N^2}{\rho^{2\alpha_N}}+\frac{t^2}{\rho^4}=1
\]
and $\norm{0}=0$.
\end{definition}

This norm is homogeneous with respect to the dilation group $(D(\lambda))_{\lambda>0}$ as long as the following 
property holds:
\[
\norm{D(\lambda)z}=\lambda\norm{z} \ \ \ \forall z\in\rnn \setminus\{0\}  \ \text{and} \ \lambda>0.
\]
Moreover the following quasi-triangular inequality holds. There exists a constant $c_T \geq 1$ such that 
\begin{equation} \label{eq-quasi-norm}
  \norm{z\circ\zeta}\leq c_T(\norm{z}+\norm{\zeta}), \qquad \norm{z^{-1}}\leq c_T \norm{z},
\end{equation}
for every $z, \zeta \in \rnn$. We denote by $d(\zeta, z) := \norm{z^{-1} \circ \zeta}$ the \emph{quasi-distance} of 
$z$ and $\zeta$. We denote by $B_\varrho (z)$  the open ball of radius $\varrho$ and center $z$ with respect to the 
above quasi-distance
\begin{equation} \label{eq-brho}
  B_\varrho (z) := \big\{ \zeta \in \rnn \mid \norm{z^{-1} \circ \zeta} < \varrho \big\}.
\end{equation}
Note that the topology induced on $\rnn$ from the norm introduced in Definition \ref{def-norm} is equivalent to the 
Euclidean one. 

\begin{remark} \label{rem-Leb-invariance}
\textnormal{The Lebesgue measure is invariant with respect to the group $(\rnn,\circ)$. Moreover, as long 
as $\det D(\lambda)=\lambda^{Q+2}$, where
\[
Q=m_0+3m_1+...+(2r+1)m_r
\]
we have that the following identity holds:
\[
    \textnormal{meas} (B_r(0))=r^{Q+2}\textnormal{meas}(B_1(0)),
\] 
where meas($B$) indicates the Lebesgue measure of the set $B$. For this reason, we will refer to $(Q+2)$ as the 
\textit{homogeneous dimension} of $\rnn$ with respect to the dilation \eqref{eq-dilation}. }

\textnormal{We also note that, in view of the structure of the matrix $B$ and of the definition of $E(\tau)$, we have 
that det$E(\tau) = 1$ for every $\tau$. In particular, the Jacobian determinant of the left translation $(x,t) \mapsto 
(\xi, \tau) \circ (x,t)$ agrees with $1$ for every $(\xi, \tau) \in \rnn$. The same is true for the Jacobian 
determinant of $(\xi, \tau) \mapsto (\xi, \tau)^{-1}$. As a consequence we have that}
\begin{equation} \label{eq-Leb-invariance}
\begin{split}
    & \int_{A} f(\zeta \circ z) d z =  \int_{\zeta \circ A} f(w) d w, 
    \qquad \quad \zeta \circ A := \big\{ w= \zeta \circ z \mid z \in A \big\},\\
    & \int_{A} f\big(z^{-1}\big) d z =  \int_{A^{-1}} f(w) d w, 
    \qquad \quad A^{-1} := \big\{ w= z^{-1} \mid z \in A \big\},\\
    & \int_{A} f(z^{-1} \circ \zeta) d z =  \int_{A^{-1}\circ \zeta} f(w) d w, 
    \quad A^{-1} \circ \zeta := \big\{ w= z^{-1}\circ \zeta \mid z \in A \big\}.
\end{split}
\end{equation}
\end{remark}

The following results is the analogous of \eqref{eq-pa>n} and \eqref{eq-pa<n} in the setting of the homogeneous Lie 
group $\left(\rnn,\circ, (D(\lambda))_{\lambda>0} \right)$. 
\begin{lemma} \label{lem-ABC} 
Let $K_{\alpha}$ denote any continuous  function which is homogeneous of degree $- \alpha$  with respect to 
$(D(\lambda))_{\lambda>0}$, for some $\alpha$ such that  $0 < \alpha < Q + 2$. For every $R>0$ we have that:
  \begin{itemize}
    \item $K_\alpha \in L^q\big(\big\{ z \in \rnn \mid \|z\| \le R \big\}\big)$ if, and only if, $q > 
\dfrac{Q+2}{\alpha}$. 
Moreover, there exists a positive constant $\widetilde c_{\alpha,q}$, only depending on $K_\alpha, 
(D(\lambda))_{\lambda>0}$ and $q$, such that 
\begin{equation}\label{eq-pa>Q+2}
  \left\| K_\alpha \right\|_{L^q(\{ z \in \rnn \mid \|z\| \le R\})} \le 
  \widetilde c_{\alpha,q} R^{\frac{Q+2}{q} - \alpha},
\end{equation}
    \item $K_\alpha \in L^q\big(\big\{ z \in \rnn \mid \|z\| \ge R \big\}\big)$ if, and only if, $q < 
\dfrac{Q+2}{\alpha}$. 
Moreover, there exists a positive constant $\widetilde c_{\alpha,q}$, only depending on $K_\alpha, 
(D(\lambda))_{\lambda>0}$ and $q$, such that 
\begin{equation}\label{eq-pa<Q+2}
  \left\| K_\alpha \right\|_{L^q(\{ z \in \rnn \mid \|z\| \ge R \})} \le 
  \widetilde c_{\alpha,q} R^{\frac{Q+2}{q} - \alpha}.
\end{equation}
  \end{itemize}
\end{lemma}

\medskip

\noindent \textit{Proof.} We compute the integrals by using the ``polar coordinates'' 
\begin{equation*} %\label{Qpolari}
\begin{cases}
x_1=\rho^{\alpha_1}\cos{\psi_1}...\cos{\psi_{N-1}}\cos{\psi_N} \\
x_2=\rho^{\alpha_2}\cos{\psi_1}...\cos{\psi_{N-1}}\sin{\psi_N} \\
 \vdots \\
x_N=\rho^{\alpha_N}\cos{\psi_1}\sin{\psi_2} \\
t=\rho^2\sin{\psi_1}.
\end{cases}
\end{equation*}
Note that, in accordance with the Definition \ref{def-norm}, the Jacobian determinant of the above change of coordinate 
is homogeneous of degree $Q+1$ with respect to the variable $\rho$, that is  $J(\rho, \psi) = \rho^{Q+1} J(1, \psi)$. 
The claim then follows by proceedings as in the Euclidean case. \hfill $\square$

\subsection{Preliminary results on convolutions in homogeneous Lie groups}
% \subsection{Continuous and compact embeddings for convolutions with homogeneous kernels}

We recall some facts concerning the convolution of functions in homogeneous Lie groups. We refer to the work of Folland 
\cite{dieci}, and to its bibliography, for a comprehensive treatment of this subject. The first result is a Young 
inequality in the non-Euclidean setting 
\begin{theorem} \label{Qyoung} Let $p, q, r \in [1,+\infty]$ be such that:
\begin{equation} \label{eq-pqr}
  1 + \frac{1}{q}=\frac{1}{p}+\frac{1}{r}.
\end{equation}
If $f\in\sobb{p}{\rnn}$ and $g\in\sobb{r}{\rnn}$, then the function $f* g$ defined as:
\[
f* g (z) := \integrale{\rnn}{}{f(\zeta^{-1}\circ z)g(\zeta)}{\zeta}
\]
belongs to $\sobb{q}{\rnn}$ and it holds:
\[
\norm{f* g}_{L^q(\rnn)} \leq \norm{f}_{L^p(\rnn)}\norm{g}_{L^r(\rnn)}.
\]
\end{theorem}

The following two theorems are the counterpart of Theorem \ref{Th-weak-Young} and Lemma \ref{14} in Section 2, 
respectively.

\begin{theorem} \label{continuitanon} {\rm (Proposition (1.11) in \cite{dieci})}
Let $K_{\alpha}$ be a continuous function, homogeneous of degree $-\alpha$ with $0 < \alpha < Q+2$,  with respect to 
the dilation \eqref{eq-dilation}. Then, for every $p \in ]1, + \infty[$, the convolution $u$ of $K_{\alpha}$ with a 
function $f \in L^p(\rnn)$ 
\begin{equation} \label{eq-convolution-sec3}
  u(z) = \int_{\rnn} K_{\alpha}(\zeta^{-1} \circ z) f(\zeta) \, {\rm d} \zeta,
\end{equation}
is defined for almost every $z \in \rnn$ and is a measurable function. Moreover
% \[
% u(z)=\integrale{\rnn}{}{\bmicr K_{\alpha} \emicr (\zeta^{-1}\circ z)f(\zeta)}{\zeta}
% \]
%Then, for every $p \in ]1, + \infty[$ 
there exists a constant $\widehat{C}_p = \widehat{C}_p(p, Q)$ such that
\[
\|u\|_{L^q(\rnn)} \leq \widehat{C}_p 
\max_{\norm{z}=1}{\abs{K_{\alpha}(z)}}\norm{f}_{L^p(\rnn)},
\]
for every $f\in\sobb{p}{\rnn}$, where  $q$ is defined by 
\[
1 + \frac{1}{q}=\frac{1}{p} + \frac{\alpha}{Q+2}.
\]
\end{theorem}

For the proof of the next result we refer to Proposition (1.11) in \cite{dieci} or Lemma 5.1 in \cite{otto}.
% It gives an estimate for homogeneous functions with respect to the dilation group. It will be needed in the 
% sequel, when dealing with uniform continuity in mean of convolutions with homogeneous kernels.

\begin{theorem} \label{Qlagrange}
Let $K_\alpha \in C^1(\rnn\setminus\{0\})$ be a homogeneous function of degree $-\alpha$ with respect to 
the group $(D(\lambda))_{\lambda>0}$. Then there exist two constants $\kappa > 1$ and $M_\alpha >0$  such that:
\[
\vert K_\alpha(\zeta)-K_\alpha(z)\vert\leq M_\alpha \, \frac{\norm{z^{-1}\circ\zeta}}{\norm{z}^{\alpha+1}}
\]
for all $z,\zeta$ such that $\norm{z}\geq \kappa \norm{z^{-1}\circ \zeta}$.
\end{theorem}

% 
% \subsection{Continuity estimates for convolutions with homogeneous kernels}
% \bmicr The continuous embedding has been proved in \cite{nove}. We recall it here for the sake of completeness. \emicr

\subsection{Compactness estimates for convolutions with homogeneous kernels}

\begin{theorem} \label{uniformecontinuitanon}
Let $K_{\alpha}$ be a $C^1(\rnn\setminus{\{0\}})$ homogeneous function of degree $-\alpha$ with $0 < \alpha 
< Q+2$, with respect to the dilation \eqref{eq-dilation}. 
% , and let $c_\alpha$ be  
% \[
% \abs{K(z)}\leq c_{\alpha}\norm{{z}}^{-\alpha}, \ \ \ c_{\alpha}>0.
% \]
Then for every $p, q \ge 1$ such that $q > p$ and 
\begin{equation} \label{eq-pqalpha-bis}
1 - \frac{\alpha + 1}{Q + 2} < \frac{1}{p} - \frac{1}{q} < 1 - \frac{\alpha}{Q + 2},
\end{equation}
there exists a positive constant $\widetilde{C}_{p,q}$ depending on $\alpha, p, q$ and on the dilation group 
$\left( D(\lambda) \right)_{\lambda >0}$ such that 
\[
\norm{u(\cdot \circ h)-u}_{L^q(\rnn)} \leq \widetilde{C}_{p,q} \norm{h}^{\frac{Q + 2}{r} - \alpha} 
\norm{f}_{L^p(\rnn)},
\]
for every $f\in \sobb{p}{{\rnn}}$, and $h \in \rnn$. Here $r$ is the constant defined by \eqref{eq-pqr}, that is
$$
  \displaystyle  1 + \frac{1}{q} = \frac{1}{r} + \frac{1}{p},
$$
and the exponent $\frac{Q + 2}{r} - \alpha$ is strictly positive, because of \eqref{eq-pqalpha-bis}. 
% introduced in \eqref{esponente-r} 
\end{theorem}

\noindent \textit{Proof.} We proceed as in the proof of Theorem \ref{Uniform continuity-bis}. We choose $z, h \in \rnn$ 
and we let 
\begin{equation} \label{eq-vz}
  v(z) := u(z \circ h) - u(z).
\end{equation}
By the formula \eqref{eq-convolution-sec3} we have $v(z) = \widetilde{I}_A(z) + \widetilde{I}_B(z) + 
\widetilde{I}_C(z)$, where
\begin{eqnarray}
 \widetilde{I}_A(z) &= & \integrale{\lbrace \zeta\in\rnn:\|\zeta^{-1} \circ z \circ h \| \geq 
\kappa \|h\|\rbrace}{}{\left( K_{\alpha}(\zeta^{-1} \circ z \circ h)-K_{\alpha}(\zeta^{-1} \circ 
z)  \right){f(\zeta)}}{\zeta}, \nonumber \\ 
 \widetilde{I}_B(z) &= & \integrale{\lbrace \zeta\in\rnn:\|\zeta^{-1} \circ z \circ h \| < 
\kappa \|h\|\rbrace}{}{{K_{\alpha}(\zeta^{-1} \circ z \circ h)}{f(\zeta)}}{\zeta}, 
\label{usplittildeQ}\\
 \widetilde{I}_C(z) &= & \integrale{\lbrace \zeta\in\rnn:\|\zeta^{-1} \circ z \circ h \| < 
\kappa \|h\|\rbrace}{}{{- K_{\alpha}(\zeta^{-1} \circ z)}{f(\zeta)}}{\zeta}. \nonumber
% \widetilde{I}_A(z) + \widetilde{I}_B(z) + \widetilde{I}_C(z)
\end{eqnarray}
Then, as in the proof of Theorem \ref{Uniform continuity-bis}, we find 
\begin{equation*}
\begin{split}
 | \widetilde{I}_A (z) | & \le c_T M_\alpha \|h\|   
 \integrale{\lbrace \zeta\in\rnn:\|\zeta^{-1} \circ z \circ h \| \geq 
\kappa \|h\|\rbrace}{}{\frac{1}{\|\zeta^{-1} \circ z \circ h\|^{\alpha + 1}}\abs{f(\zeta)}}{\zeta}, \\
 | \widetilde{I}_B (z) | & \le c_\alpha \integrale{\lbrace \zeta\in\rnn:\|\zeta^{-1} \circ 
z \circ h \| < \kappa \|h\|\rbrace}{}{\frac{1}{\|\zeta^{-1} \circ z \circ h\|^{\alpha}}\abs{f(\zeta)}}{\zeta}, \\
 | \widetilde{I}_C (z) | & \le c_\alpha \integrale{\lbrace \zeta\in\rnn:\|\zeta^{-1} \circ z 
\circ h \| < \kappa \|h\|\rbrace}{}{\frac{1}{\|\zeta^{-1} \circ z\|^{\alpha}}\abs{f(\zeta)}}{\zeta}. 
\end{split}
\end{equation*}
The first estimate follows from Theorem \ref{Qlagrange} and the constant $c_T$ is the one appearing in 
\eqref{eq-quasi-norm}, while $c_\alpha := \max_{\norm{w}=1} K_\alpha(w)$. 

We next compute the $L^r$ norm of the homogeneous functions appearing above. In view of Remark 
\ref{rem-Leb-invariance} and Lemma \ref{lem-ABC}, we have that
\begin{equation*}
%\begin{split}
\left\| \frac{1}{\norm{\zeta^{-1} \circ z \circ h}^{\alpha + 1}} \right\|_{L^r(\lbrace 
\zeta\in\rnn:\|\zeta^{-1} \circ z \circ h \| \geq 
\kappa \|h\|\rbrace)} = \widetilde{C}_A(r, Q) \|h\|^{\frac{Q+2}{r} - (\alpha + 1)},
%\end{split}
\end{equation*}
where $r$ is the exponent introduced in \eqref{esponente-r} and $\widetilde{C}_A(r,Q)$ is a constant depending on $Q$ 
and $r$  (hence on $Q, p, q$), and on the dilation group $\left( D(\lambda) \right)_{\lambda >0}$. 
Using again Remark \ref{rem-Leb-invariance} and Lemma \ref{lem-ABC}, and  we also find 
\begin{equation*}
%\begin{split}
\left\| \frac{1}{\norm{\zeta^{-1} \circ z \circ h}^{\alpha}} \right\|_{L^r(\lbrace 
\zeta\in\rnn:\|\zeta^{-1} \circ z \circ h \| \le 
\kappa \|h\|\rbrace)} = \widetilde{C}_B(r, Q) \|h\|^{\frac{Q+2}{r} - \alpha},
%\end{split}
\end{equation*}
The same argument applies to $\widetilde{I}_C$, by using the quasi-triangular inequality \eqref{eq-quasi-norm}, so that 
\begin{equation*}
%\begin{split}
\left\| \frac{1}{\norm{\zeta^{-1} \circ z \circ h}^{\alpha}} \right\|_{L^r(\lbrace 
\zeta\in\rnn:\|\zeta^{-1} \circ z \circ h \| \le 
\kappa \|h\|\rbrace)} = \widetilde{C}_C(r, Q) \|h\|^{\frac{Q+2}{r} - \alpha},
%\end{split}
\end{equation*}
Note that the three above integrals converge if, and only if,
\begin{equation} \label{cond-su-r-Q}
\frac{Q+2}{\alpha + 1} < r < \frac{Q+2}{\alpha}.
\end{equation}
We note that \eqref{eq-pqalpha-bis} is equivalent to \eqref{cond-su-r-Q} and that the second inequality in  
\eqref{cond-su-r-Q} says that the exponent $\frac{Q + 2}{r} - \alpha$ appearing in the statement of this Theorem is 
strictly positive. By the Young inequality (Theorem \ref{Qyoung}) we conclude that
\begin{equation*}
\begin{split}
& \left\| \widetilde{I}_A \right\|_{L^q(\rnn)} \le c_T M_\alpha \widetilde{C}_A(r,Q) \|h\|^{\frac{Q+2}{r}- \alpha} 
\left\| f \right\|_{L^p(\rnn)}, \\
& \left\| \widetilde{I}_B \right\|_{L^q(\rnn)} \le c_\alpha \widetilde{C}_B(r,Q) \|h\|^{\frac{Q+2}{r} - 
\alpha} \left\|f\right\|_{L^p(\rnn)}, \\
& \left\| \widetilde{I}_B \right\|_{L^q(\rnn)} \le c_\alpha \widetilde{C}_C(r,Q) \|h\|^{\frac{Q+2}{r} - 
\alpha} \left\|f\right\|_{L^p(\rnn)}.
\end{split}
\end{equation*}
The thesis is obtained with $\widetilde{C}_{p,q} := c_T M_\alpha \widetilde{C}_A(r,Q) +  c_\alpha \widetilde{C}_B(r,Q) 
+ c_\alpha \widetilde{C}_C(r,Q)$. \hfill $\Box$

\begin{remark} \label{liston-b}
Similarly as we did in Remark \ref{listona}, we can state the conditions on $q$ for the validity of 
\eqref{eq-pqalpha-bis}. They can can be obtained by substituting the \emph{dimension} $n$ with the 
\emph{homogeneous dimension} $Q+2$. We explicitly write here the condition for $\alpha = Q$ and $\alpha = Q+1$, which 
occur in the representation formulas for the solutions to $\L u = f$. Moreover, when $\alpha = Q+1$, we only consider 
the case $1 < p < Q+2$, as we apply Theorem \ref{uniformecontinuitanon} to prove the compactness of the embedding of
Theorem \ref{continuitanon}, which holds only for $p < Q+2$. For the same reason, when $\alpha = Q$, we only consider 
the case $\displaystyle 1 < p < \frac{Q + 2}{2}$.  
\begin{itemize}
\item[(i)]  If {$\alpha = Q$} and $\displaystyle 1 < p < \frac{Q + 2}{2}$, we have that $p^* < q < p^{**}$,
% If \fbox{$\alpha = Q$}
% \begin{itemize}
% \item[(i.1)] If $\displaystyle 1 < p < \frac{Q + 2}{2}$ then 
% \[  
% \frac{p (Q + 2)}{(Q + 2) - p} < q < \frac{p (Q + 2)}{(Q + 2) - 2p}.
% \]
% \item[(i.2)] If $\displaystyle \frac{Q + 2}{2} \le p \le Q + 2$ then 
% \[ 
% q > %\max \left \{1, \frac{p (Q + 2)}{(Q + 2) - p} \right \} = 
% \frac{p (Q + 2)}{(Q + 2) - p}.
% \]
% \end{itemize}
\item[(ii)] If {$\alpha = Q + 1$} and $1 < p < Q+2$ we have that $p < q < p^*$.
% \begin{itemize}
% \item[(ii.1)] If $\displaystyle 1 \le p \le Q + 2$ then 
% \[  
% p < q < p^* = \frac{p (Q + 2)}{(Q + 2) - p}.
% \]
% \item[(ii.2)] If $p \ge Q + 2$ then 
% \[ 
% q > p.
% \]
% \end{itemize}
\end{itemize}
\end{remark}

\subsection{Proof of our main result}

\noindent {\it Proof of Theorem \ref{main-th-1}}. The proof of the Sobolev inequality is a direct consequence of 
Theorem \ref{continuitanon}, with $\alpha = Q+1$ when considering $v_1, \dots, v_{m_0}$, and $\alpha = Q$ as we 
consider $u$. 

The compactness of the embedding is a direct consequence of Theorem \ref{uniformecontinuitanon}. As noticed in Remark 
\ref{liston-b}, it applies to the derivatives $\partial_{\xi_j}\Gamma$, that are homogeneous functions of degree 
$-(Q+1)$, only when $p < q < p^*$. Moreover, a direct computation based on \eqref{eq-pqr} shows that
\begin{equation}
  \frac{Q + 2}{r} - (Q+1) =  (Q+2) \left(\frac{1}{q}- \frac{1}{p^*} \right).
\end{equation}
Analogously, as $\Gamma$ is homogeneous of degree $-Q$, we need to consider $p^* < q < p^{**}$. In this case, by using 
again \eqref{eq-pqr} we find
\begin{equation}
  \frac{Q + 2}{r} - Q =  (Q+2) \left(\frac{1}{q}- \frac{1}{p^{**}} \right).
\end{equation}

The proof of the Morrey embedding is obtained by the same argument used in the proof of Theorem \ref{e}. 
We consider the  function $v(z)= u(z \circ h) - u(z)$ introduced in  \eqref{eq-vz} and, as in the proof of  Theorem 
\ref{uniformecontinuitanon}, we write $v(z) = \widetilde{I}_A(z) + \widetilde{I}_B(z) + \widetilde{I}_C(z)$, where the 
functions $\widetilde{I}_A, \widetilde{I}_B, \widetilde{I}_C$ are defined in \eqref{usplittildeQ}. The conclusion of 
the 
proof is obtained by using the Young inequality stated in Theorem \ref{Qyoung}. \hfill $\Box$

\section{Representation formulas}
\setcounter{section}{4} \setcounter{equation}{0}
\setcounter{theorem}{0} 

\subsection{Fundamental solution to $\L_{0}$ and representation formula}

In this Section we focus on the representation formulas for the equation $\L u = 0$, for $\L$ satisfying  the 
assumptions (H.1) and (H.2), then we prove Theorem \ref{main-th-2}.

We first recall the definition of \emph{weak solution} to $\L u = 0$, then we recall that, 
under these assumptions, the fundamental solution to $\L_0$ has been derived by H\"ormander in \cite{dieci}. 
We say that $u$ is a weak solution to $\L u = 0$ in an open set $\Omega\subset\rnn$ if 
$u, \partial_{x_1}u, \dots, \partial_{x_{m_0}}u, Yu \in L_{\mathrm{loc}}^2(\Omega)$ and 
\[
\integrale{\Omega}{}{-\langle A(z)\nabla u(z),\nabla\psi(z)\rangle +\psi(z) Y u(z)}{z}=0,
\]
for all $\psi\in C_0^\infty(\Omega)$. 

% \begin{definition}
% We say that $u$ is a weak solution to $\L u = 0$ in an open set $\Omega\subset\rnn$ if 
% $u, \partial_{x_1}u, \dots, \partial_{x_{m_0}}u, Yu \in L_{\mathrm{loc}}^2(\Omega)$ and 
% \[
% \integrale{\Omega}{}{-\langle A(z)\nabla u(z),\nabla\psi(z)\rangle +\psi(z) Y u(z)}{z}=0,
% \]
% for all $\psi\in C_0^\infty(\Omega)$.
% \end{definition}
With the notations introduced in Section 3, we let
% An explicit estimate of the fundamental solution to $\L_0$ has been derived in \cite{dieci}. Under
% \begin{theorem} \label{th-fund-sol}
\[
C(t)=\integrale{0}{t}{E(s)A_0 E^T(s)}{s}
\]
The assumptions (H.1) and (H.2) guarantee that the matrix $C(t)$ is strictly positive for all $t>0$. In this case its 
inverse $C^{-1}(t)$ is well defined and the fundamental solution to $\L_{0}$ with singularity at the origin of $\rnn$, 
is given by 
\begin{equation} \label{eq-fund-sol-K}
  \Gamma((x,t),(0,0))=
\begin{cases}
\dfrac{(4\pi)^{-\frac{N}{2}}}{\sqrt{\mathrm{det}C(t)}}\mathrm{exp}
\left(-\frac{1}{4}\left\langle C^{-1}(t)x,x\right\rangle \right) , & \mbox{if } t>0, \\ 
0, & \mbox{if } t\leq 0.
\end{cases}
\end{equation}
To simplify the notation, in the sequel we will write $\Gamma(x,t,\xi, \tau)$ instead of 
$\Gamma((x,t),(\xi,\tau))$, and $\Gamma(x,t)$ instead of $\Gamma(x,t,0,0)$.
The fundamental solution $\Gamma(x,t,\xi,\tau)$ of $\L_{0}$ with pole at $(\xi,\tau)$ is the ``left 
translation'' of $\Gamma(\cdot,0,0)$ with respect to the group $(\mathbb{R}^{N+1},\circ)$:
\[
\Gamma(x,t,\xi,\tau)=\Gamma((\xi,\tau)^{-1}\circ(x,t),0,0).
\]
%\end{theorem}
Let us explicitly remark that $\Gamma(\cdot,0,0)$ is homogeneous of degree $-Q$ with respect to the group 
$(D(\lambda))_{\lambda>0}$ and $\partial_{x_j}\Gamma(\cdot,0,0)$ is homogeneous of degree $-Q-1$, for $j=1,...,m_0$.
Moreover, also $\partial_{\xi_j}\Gamma(0,0, \cdot)$ is homogeneous of degree $-Q-1$, for $j=1,...,m_0$.

We next represent weak solutions to $\L u = 0$ as convolutions with the fundamental solution $\Gamma$ to $\L_0$ and 
to its derivatives $\partial_{\xi_1}\Gamma, \dots, \partial_{\xi_{m_0}}\Gamma$.

% and $r,s\in\mathbb{R}$, with $0<s<r$. Let $\varphi\in C^\infty(\mathbb{R})$ be 
% a cut-off function defined as:
% \[
% \varphi(t)=1 \ \mathrm{for \ all} \ 0\leq t\leq s, \ \ \ \varphi(t)=0 \ \mathrm{for \ all} \ t\geq r
% \]
% For all $\zeta_0\in\Omega$ and $r>0$ such that $B_r(\zeta_0)\subset\Omega$ we set:
% \begin{equation} \label{eta}
% \eta(z)=\varphi(\norm{\zeta_0^{-1}\circ z})
% \end{equation}
Consider any open set $\Omega \subseteq \rnn$ let $u$ be a function such that $u, \partial_{x_1} u, \dots, 
\partial_{x_{m_0}} u, Yu  \in L^p_{\text{\rm loc}}(\Omega)$, and $\eta \in C_0^\infty(\Omega)$ is any cut-off 
function, then, by an elementary density argument we find
\begin{equation*} %\label{rappresentazione_L0}
\begin{split}
(\eta u)(z) = & -\integrale{\mathbb{R}^{N+1}}{}{[\Gamma(z,\cdot)\L_0 (\eta u)](\zeta)}{\zeta} = \\
  & \integrale{\mathbb{R}^{N+1}}{}{[\langle A_0 \nabla_\xi \Gamma(z,\cdot),  \nabla (\eta u) 
  \rangle](\zeta)}{\zeta}-\integrale{\mathbb{R}^{N+1}}{}{[\Gamma(z,\cdot)Y(\eta u)](\zeta)}{\zeta}.
\end{split}
\end{equation*}
If moreover $u$ is a weak solutions to $\L u = \text{div} (A_0 F) + f$, then $\L_0 u =$ div$\left((A_0 - A\right) 
\nabla u + A_0 F) + f$, for some $f \in L^p_\text{\rm loc}(\Omega)$ and some vector valued function $F = (F_1, \dots, 
F_{m_0}, 0, \dots, 0)$ with $F_1, \dots, F_{m_0} \in L^p_\text{\rm loc}(\Omega)$, then we obtain the following 
representation formula introduced in Theorem 3.1 of \cite{otto}, and used in the proof of Theorem 3.3 of \cite{sei}. 
\begin{theorem} \label{th-rappr} 
If $u$ is a weak solution to $\L u = \text{ \rm div} (A_0 F) + f $, 
%and $\eta$ is the function introduced in \eqref{eta}, 
in some open set $\Omega \subset \rnn$, with $f, F_1, \dots, F_{m_0} \in L^p_\text{\rm loc}(\Omega)$, and $\eta$ is the 
cut-off function defined above,  then:
\begin{equation} \label{rappresentazione}
\begin{split}
(\eta u)(z) = & \integrale{\mathbb{R}^{N+1}}{}{ \! \!  [ \eta \langle \nabla_\xi \Gamma(z,\cdot), 
(A_0 - A) \nabla u + A_0 F \rangle](\zeta)}{\zeta} -
\integrale{\mathbb{R}^{N+1}}{}{\! \! [\Gamma(z,\cdot) (\langle A \nabla \eta , \nabla u \rangle + \eta f )]}{\zeta} \\
 + & \integrale{\mathbb{R}^{N+1}}{}{[\langle A_0 \nabla_\xi \Gamma(z,\cdot),  \nabla \eta \rangle u ](\zeta)}{\zeta}-
\integrale{\mathbb{R}^{N+1}}{}{\Gamma(z,\cdot)(Y \eta) u}{\zeta}.
\end{split}
\end{equation}
\end{theorem} 

In the following statement $B_{\varrho}(z_0)$ denotes the ball defined in \eqref{eq-brho}, and $c_T$ is the constant in 
\eqref{eq-quasi-norm}.

\begin{proposition} \label{prop-estimates} Let $\Omega$ be an open set of $\rnn$, and let $u$ be a weak solution to  
$\L 
u = \text{div} (A_0 F) + f$ in $\Omega$. Suppose that $u, f, \partial_{x_1} u, \dots, \partial_{x_{m_0}} u, F_1, 
\dots, F_{m_0} \in L^p_\text{\rm loc}(\Omega)$. Then for every $z_0 \in \Omega$, and $\varrho, \sigma >0$ such that the 
ball $B_{\varrho}(z_0)$ is contained in $\Omega$, and $\sigma < \dfrac{\varrho}{2 c_T}$, we have:
\begin{itemize}
  \item (Sobolev embedding) if $1 < p < Q+2$, then there exists a positive constant $C_p$ such that
\begin{equation*} %\label{eq-Young-main}
\begin{split}
  \| u \|_{L^{p^*}(B_{\sigma}(z_0))} \le \, &  C_{p} \big( \| u \|_{L^{p}(B_{\varrho}(z_0))} + 
  \| f \|_{L^{p}(B_{\varrho}(z_0))} \\ 
  + & \, \| A_0 \nabla u \|_{L^{p}(B_{\varrho}(z_0))} + \|  A_0 F \|_{L^{p}(B_{\varrho}(z_0))} \big), 
\end{split}
\end{equation*}
 \item (Compactness) if moreover $p < q < p^*$, then there exists a positive constant 
$\widetilde C_{p, q}$ such that 
\begin{equation*} %\label{eq-Sobolev-unif-main}
\begin{split}
  \| u ( \cdot \circ h) - u \|_{L^{q}(B_{\sigma}(z_0))} \le \, & \widetilde C_{p, q}  
  \big( \| u \|_{L^{p}(B_{\varrho}(z_0))} +  \| f \|_{L^{p}(B_{\varrho}(z_0))} \\
 + & \, \| A_0 \nabla u \|_{L^{p}(B_{\varrho}(z_0))} + \|  A_0 F \|_{L^{p}(B_{\varrho}(z_0))} \big) \
  \|h\|^{(Q+2) \left( \frac{1}{q}- \frac{1}{p^*} \right)}, 
\end{split}
\end{equation*}
for every $ h \in B_{\sigma}(z_0)$, 
\item (Morrey embedding) if  $p > Q+2$, then there exists a positive constant $\widetilde C_p$ such that
\begin{equation*} %\label{eq-Morrey-main}
\begin{split}
   | u (z) - u (\zeta) | \le \, & \widetilde C_{p}  \big( \| u \|_{L^{p}(B_{\varrho}(z_0))} + 
   \| f \|_{L^{p}(B_{\varrho}(z_0))} \\   
   + & \,  \| A_0 \nabla u \|_{L^{p}(B_{\varrho}(z_0))} + \|  A_0 F \|_{L^{p}(B_{\varrho}(z_0))} \big) 
   \|\zeta^{-1} \circ z\|^{1 - \frac{Q+2}{p}}, 
\end{split}
\end{equation*}
for every $\ z, \zeta \in B_{\sigma}(z_0)$.
 \end{itemize}
\end{proposition}

\medskip

\noindent {\it Proof}. We apply Theorem \ref{th-rappr} with a function $\eta$ supported in the ball $B_{\varrho}(z_0)$ 
and such that $\psi(z) = 1$ for every $z \in B_{2 c_T \sigma}(z_0)$. It is not difficult to check that a cut-off 
function with the above properties exists (see formula (3.3) in \cite{otto} for instance). Note that the integrals 
appearing in the equation \eqref{rappresentazione} involving $\partial_{x_1}u, \dots, \partial_{x_{m_0}}u, F_1, \dots, 
F_{m_0}$ are convolutions of $\partial_{\xi_1}\Gamma, \dots, \partial_{\xi_{m_0}}\Gamma$, that are 
homogeneous kernels of degree $-(Q+1)$, with functions belonging to $L^p(B_{\varrho}(z_0))$ multiplied by bounded 
functions compactly supported in $B_{\varrho}(z_0)$.  For these terms of the representation formula  
\eqref{rappresentazione} the thesis then follows from a direct application of Theorem \ref{main-th-1}. Indeed, by our 
choice of $\eta$, we have $u(z) = (\eta u)(z)$ for every $z \in B_{\sigma}(z_0)$. Moreover, if $z, h \in 
B_{\sigma}(z_0)$ we also have that $z \circ h \in B_{2 c_T \sigma}(z_0)$, by \eqref{eq-quasi-norm}, then also $u(z 
\circ h) = (\eta u)(z \circ h)$.

We next consider the terms involving $u$ and $f$. They are convolutions of $\Gamma$, that is a 
homogeneous kernel of degree $-Q$, with $u$ and $f$, multiplied by bounded functions compactly supported in 
$B_{\varrho}(z_0)$. Moreover, $u$ and $f$ belong to $L^r(B_{\varrho}(z_0))$, for every $r$ such that $1 \le r \le p$. 
We then choose $r$ such that 
\begin{equation*}
  \frac{1}{r} = \frac{1}{p} + \frac{1}{Q+2}
\end{equation*}
and we apply again Theorem \ref{main-th-1} with $p$ replaced by $r$. This concludes the proof. \hfill $\Box$

\medskip
\noindent {\it Proof of Theorem \ref{main-th-2}}. It follows from Proposition \ref{prop-estimates} by a simple covering 
argument. The constant $\widetilde \varrho$ can be chosen as follows. We let 
$$
\overline \varrho : = \min \big\{\varrho >0 \mid B_\varrho(z) \subset\Omega \ \text{for every} \ z \in \K \big\},
$$
then $\widetilde \varrho := \dfrac{\overline \varrho}{3 c_T}$, so that we can choose $\sigma = \widetilde \varrho$ in 
every ball of the covering of $\K$. \hfill $\Box$

\medskip
\noindent {\it Proof of Theorem \ref{main-th-3}}. If $f$ is a weak solution to \eqref{eq-kinetic}, then it is a weak 
solution to 
\begin{equation} \label{eq-kinetic-2}
  \partial_t f + \langle v, \nabla_x f \rangle = \Delta_v f + g - \text{\rm div}_v G, 
\end{equation}
where $G_j = \partial_{v_j} f, j=1, \dots, n$. Note that the homogeneous dimension of the operator in 
\eqref{eq-kinetic-2} is $Q+2 = n + 3n +2$. By our assumptions $G_j \in L^p (\Omega)$ for every $j=1, \dots, n$. 
Then the proof can be concluded by the same argument used in the proof of Theorem \ref{main-th-2}. \hfill $\Box$

\section{Conclusion}
\setcounter{section}{5} \setcounter{equation}{0}
\setcounter{theorem}{0} 

The method used in this article for Kolmogorov equations can be adapted to the study of a wider family of differential 
operators, provided that they have a fundamental solution and that are invariant with respect to a suitable Lie group 
on their domain. Sobolev inequalities for operators of this kind have been proven in \cite{quindici}. We recall here 
the assumptions on the operators.

Consider a differential operator in the form
\begin{equation}\label{e1-parabolic}
     \L u : = \sum_{i,j=1}^m X_j \left( a_{ij}(x,t) X_i u \right) + X_0 u - \p_t u,
\end{equation}
where $(x,t)=(x_{1},\dots,x_{N},t)$ denotes the point in $\R^{N+1}$, and $1 \le m \le N$. The $X_{j}$'s
in \eqref{e1-parabolic} are smooth vector fields acting on $\R^{N}$, i.e.
\begin{equation*}
     X_{j}(x,t)=\sum_{k=1}^{N} b_{k}^{j}(x,t)\p_{x_{k}}, \qquad j=0,\dots,m,
\end{equation*}
and every $b_{k}^{j}$ is a $C^{\infty}$ function. In the sequel we always denote by $z=(x,t)$ the point
in $\rnn$, and by $A$ the $m \times m$ matrix $A = \left( a_{i,j} \right)_{i,j = 1, \dots, m}$. We also consider the 
\emph{elliptic} analogous of $\L$
\begin{equation}\label{e1-elliptic}
     \L u : = \sum_{i,j=1}^m X_j \left( a_{ij}(x) X_i u \right)
\end{equation}
In both cases we assume that the coefficients of the matrix $A$ are bounded measurable functions, and that $A$ is 
symmetric and uniformly positive, that is, there is a positive constant $\mu$ such that
\begin{equation*} %\label{e1-unif-elliptic}
      \sum_{i,j=1}^m a_{ij}(x,t) \xi_i, \xi_j \ge \mu |\xi|^2, \qquad \text{for every} \ \xi \in \R^m, 
\end{equation*}
and for every $(x,t) \in \rnn$ (or for every $x \in \rn$ as we consider the operator $\L$ in \eqref{e1-elliptic}). 

Clearly, the Laplace operator $\Delta$ and the heat operator $\Delta - \partial_t$ write in the form 
\eqref{e1-elliptic} and \eqref{e1-parabolic}, respectively, if we choose $X_j := \partial_{x_j}$ for $j=1, \dots, N$,  
$X_0 := 0$, and the matrix $A$ agrees with the $N \times N$ identity $I_N$. 
In the sequel we will use the following notations:
\begin{equation*} %\label{eY}
  X = \left( X_{1}, \dots, X_{m} \right), \qquad Y=X_{0}-\p_{t},
\qquad  \div F = \sum_{j=1}^m X_j F_j,
\end{equation*}
for every vector field $F = (F_{1}, \dots, F_{m} )$, so that the expression $\L u$ reads
\begin{equation*}
  \L u = \div (A Xu) + Y u.
\end{equation*}
Finally, when $A$ is the $m \times m$ identity matrix, we will use the notation
\begin{equation*} %\label{e0}
  \LO:= \sum_{k=1}^{m}X_{k}^{2}+Y.
\end{equation*}

%\appendix
%\input{Bibliog}

%\bibliography{bib}

\begin{thebibliography}{10} 
%\addcontentsline{toc}{chapter}{Bibliografia}
\bibitem{tre} R. A. Adams e J. J. F. Fournier, \textit{Sobolev spaces}, Academic Press, 2003.
\bibitem{quattordici} {F. Bouchut}, \textit{Hypoelliptic regularity in kinetic equations} J. Math. Pures Appl. (9), 
81(11) (2002) 1135--1159.
\bibitem{cinque} M. Bramanti, M. C. Cerutti e M. Manfredini, \textit{$L^p$ estimates for some ultraparabolic operators 
with discontinuous coefficients}, J. Math. Anal. Appl. 200(2) (1996) 332--354.
\bibitem{due} H. Brezis, \textit{Analisi funzionale. Teoria ed applicazioni}, Liguori, 1986.
\bibitem{sette} C. Cinti , A. Pascucci, S. Polidoro, \textit{Pointwise estimates for solutions to a class of 
non-homogeneous Kolmogorov equations}, Math. Ann. 340(2) (2008) 237--264.
\bibitem{quindici} C. Cinti , S. Polidoro, \textit{Pointwise local estimates and Gaussian upper bounds for a class of 
uniformly subelliptic ultraparabolic operators} J. Math. Anal. Appl. 338 (2008) 946--969. 
\bibitem{uno} L. C. Evans, \textit{Partial Differential Equations: Second Edition}, AMS, 2010.
\bibitem{nove} G. B. Folland, \textit{Subelliptic estimates and function spaces on nilpotent Lie groups}, Ark. Mat. 
13(2) (1975) 161--207.
\bibitem{tredici} {F. Golse, C. Imbert, C. Mouhot and A. F. Vasseur}, \textit{Harnack inequality for kinetic 
Fokker-Planck equations with rough coefficients and application to the Landau Equation}  (to appear on Ann. Sc. Norm. 
Super. Pisa Cl. Sci. (5), DOI Number: 10.2422/2036-2145.201702\_001 preprint, arXiv.org:1607.08068) (2017).
\bibitem{dieci} L. H\"ormander, \textit{Hypoelliptic second order differential equations}, Acta Math. 119 (1967) 
147--171.
\bibitem{dodici} {E. Lanconelli, S. Polidoro}, \textit{On a class of hypoelliptic evolution operators} 
Rend. Sem. Mat. Univ. Politec. Torino 52,1 (1994) 29--63.
\bibitem{quattro} V. Manco, G. Metafune e C. Spina, \textit{Equazioni ellittiche del secondo ordine. Parte seconda: 
teoria $L^p$}, Università di Lecce, Quaderni di matematica 4, 2005.
\bibitem{otto} M. Manfredini, S. Polidoro, \textit{Interior regularity for weak solutions of ultraparabolic equations 
in divergence form with discontinuous coefficients}, Boll. Unione Mat. Ital. Sez. B Artic. Ric. Mat. 1(8) (1998) 
651--675.
\bibitem{sei} A. Pascucci e S. Polidoro, \textit{The Moser’s iterative method for a class of ultraparabolic equations}, 
Commun. Contemp. Math. 6 (2004) 395--417.
\bibitem{otto.1} S. Polidoro, M.A. Ragusa, \textit{H\"older regularity for solutions of an ultraparabolic equations in 
divergence form}, Potential Anal. 14 (2001) 341--350.
\bibitem{undici} {E. M. Stein},  \textit{Singular integrals and differentiability properties of functions},
{Princeton Mathematical Series, No. 30}, {Princeton University Press, Princeton, N.J.}, (1970) {xiv+290}.
\end{thebibliography}

\end{document}